\documentclass[12pt,a4paper,onecolumn,reqno]{amsart}

\usepackage{amsmath,amsfonts,amsthm,amsopn,mathrsfs,amssymb,cite}
\usepackage{epsfig,verbatim}

\usepackage{booktabs}
\usepackage{epstopdf}
\usepackage{acronym}
\usepackage{enumitem}   
\usepackage{caption}
\usepackage{subcaption}
\usepackage{tabularx}
\newtheorem{thm}{Theorem}[section]
\newtheorem{remark}[thm]{Remark}
\newtheorem{defi}[thm]{Definition}
\newtheorem{lem}[thm]{Lemma}
\newtheorem{cor}[thm]{Corollary}

\renewcommand{\H}{\mathcal{H}}

\renewcommand{\leq}{\leqslant}

\renewcommand{\geq}{\geqslant}
\renewcommand{\epsilon}{\varepsilon}
\newcommand\R{\mathcal{R}}

\newcommand\W{\mathcal{W}}
\newcommand\N{\mathcal{N}}

\newcommand{\cH}{\mathcal{H}}
\newcommand{\cU}{\mathcal{U}}
\newcommand{\cT}{\mathcal{T}}

\DeclareMathOperator{\e}{e}

\DeclareMathOperator{\sinc}{sinc}
\newcommand{\bN}{\mathbb{N}}
\newcommand{\inv}{^{-1}}

\newcommand{\arro}{\,\Rightarrow \, }
\acrodef{IPRM}[IPRM]{Inverse Polynomial Reconstruction Method}
\acrodef{GEIM}[GEIM]{Generalized Empirical Interpolation Method}
\acrodef{SNR}[SNR]{signal-to-noise ratio}

\title[Sampling and Reconstruction in Distinct Subspaces]{Sampling and Reconstruction in Distinct Subspaces Using Oblique Projections}

\author{Peter Berger}
\address{Institute of Telecommunications  \\
Vienna University of Technology, \newline
Gusshausstrasse 25/389\\
A-1040 Vienna, Austria}
\email{peter.berger@nt.tuwien.ac.at}
\author{Karlheinz Gr\"ochenig}
\address{Faculty of Mathematics \\
University of Vienna \\
Oskar-Morgenstern-Platz 1\\
A-1090 Vienna, Austria}
\email{karlheinz.groechenig@univie.ac.at}
\author{Gerald Matz}
\address{Institute of Telecommunications  \\
Vienna University of Technology, \newline
Gusshausstrasse 25/389\\
A-1040 Vienna, Austria}
\email{gerald.matz@tuwien.ac.at}

\subjclass[2000]{}
\date{}
\thanks{Funding by the Austrian Science Fund (FWF) through grant NFN
  SISE (S10602) and  P26273 - N25 and by the Vienna Science and Technology Fund (WWTF)
  through project VRG12-009 and project ICT15-119.} 

\begin{document}
\thispagestyle{plain}
\pagestyle{plain}
\begin{abstract}
We study reconstruction operators on a Hilbert space that are exact on
a given reconstruction subspace. Among those the reconstruction
operator obtained by the least squares fit has the smallest operator
norm, and therefore is most stable with respect to noisy
measurements. We then  construct  the operator with the smallest possible 
quasi-optimality constant, which is the  most
stable with respect to a systematic error appearing before the
sampling process (model uncertainty). We describe how to vary continuously between the two
 reconstruction methods, so that we can trade stability 
for quasi-optimality. As  an application we study the  reconstruction of
a compactly supported function from nonuniform samples of its  Fourier
transform. 
\end{abstract}
\maketitle

\section{Introduction}\label{intro}
\subsection{The reconstruction problem}
In this paper we treat the following sampling problem.
Let $\H$ be a separable Hilbert space over $\mathbb{C}$ with inner
product $\langle \cdot,\cdot\rangle_\mathcal{H}$.
We assume that we are given linear measurements 
$(\langle f,u_j\rangle_\mathcal{H})_{j\in\mathbb{N}}$, $u_j\in\H$, 
of an unknown function $f\in \H$. We call $(u_j)_{j\in\mathbb{N}}$ the
\textit{sampling frame} and
$\mathcal{U}:=\overline{\textnormal{span}}(u_j)_{j\in\mathbb{N}}$ the
\textit{sampling space}. 
Our goal is to approximate the function $f$ by an element in the
\textit{reconstruction space}
$\cT:=\overline{\textnormal{span}}(t_k)_{k\in\mathbb{N}}$ with
$t_k\in \H$,  by a  series expansion $\tilde{f} =
\sum_{k\in\mathbb{N}} c_k t_k$ from the given measurements.
The main point is that  in general  the reconstruction space is 
distinct from the sampling space, whereas in  classical
frame theory these two spaces coincide. 

\subsection{Areas of application and related work}

This type of sampling problem arises in many concrete applications and
in the numerical modelling of infinite dimensional problems. 

(i) \emph{Sampling of bandlimited functions.}  In~\cite{gro} a
bandlimited function is approximated from finitely many, nonuniform  samples by
means of a trigonometric polynomial. In this case the sampling space
consists of the reproducing kernels $u_j (x) = \frac{\sin \pi (x- x_j)}{\pi
  (x-x_j)}, j=1, \dots , n,$ and the reconstruction vectors are
$t_k(x) = e^{2\pi i k x/(2M+1)} \chi _{[-M,M]}(x), |k| \leq M$.

(ii) \emph{Inverse Polynomial Reconstruction Method.} In this method
one tries to approximate an algebraic polynomial or an analytic
function from its Fourier samples. Thus the sampling space consists of
vectors $u_j(x) = e^{\pi i j x}\chi _{[-1,1]}(x), j=1, \dots , m$,
and the reconstruction space consists of a suitable polynomial basis,
usually the monomials $t_k(x) = x^k, k=0, \dots , n$, or the Legendre
polynomials. This method claims to efficiently mitigate the Gibbs phenomenon
~\cite{iprm1,iprm2,iprm3,iprm4},  and, indeed, the
modified inverse polynomial reconstruction method \cite{Hrycak} leads
to a numerically stable reconstruction when $m\geq n^2$.

(iii) \emph{Fourier sampling.} More generally, the goal is to
approximate a compactly supported function in some smoothness class
from its nonuniform Fourier samples $\hat{f}(\omega _j)$. Thus the sampling space
consists again of the functions $u_j(x) = e^{\pi i \omega _j x}\chi
_{[-1,1]}(x)$. The reconstruction space depends on the signal
model and on a priori information. If $f$ is smooth and belongs to a
Besov space, then the reconstruction space may be  taken to be a
wavelet subspace. The problems of Fourier sampling have motivated
Adcock and Hansen to revisit  nonuniform sampling theory and to create
the impressive and useful framework of generalized sampling~\cite{adcock_2d_wavelet,adcock,ahc1,1,Ma2017}.  

(iv) \emph{Model reduction in parametric  partial differential
  equations and the generalized empirical interpolation method.} In general  the solution manifold  
to a parametric partial differential equation is  quite complicated,
therefore  it is approximated by  finite-dimensional spaces $\cT_n$.
The \ac{GEIM} \cite{MR3051403,MR3318671} builds an interpolant in an
$n$-dimensional space $\cT_n$ based on the knowledge of $n$ physical
measurements $(\langle f,u_j\rangle_\mathcal{H})_{j=1}^n$. In
\cite{MR3341243,binev} an extension based on a least squares 
method has been proposed, where the dimension $m$ of $\cT_m$ is smaller
than the number $n$ of the measurements $(\langle
f,u_j\rangle_\mathcal{H})_{j=1}^n$. A further generalization to Banach 
spaces is contained  in \cite{ronald}.
The focus in \cite{MR3341243,binev} lies in minimizing the error
caused by the model mismatch. This is done by using a correction term
outside of the reconstruction space, which means that (in contrast to
our work) the reconstruction is allowed to be located outside of the
reconstruction space. 
This approach is optimal in the absence of measurement noise
\cite{binev}. 

In all these problems the canonical approximation or reconstruction is
by means of a least squares fit, namely 
\begin{equation}\label{wlsq}
\tilde{f} = \underset{g\in \cT}{\textnormal{arg~min}}
\sum_{j\in\mathbb{N}} w_j\left|\langle g, u_j\rangle_\mathcal{H}-
  d_j\right|^2 \, .
\end{equation}
The weights   $w_j$ are usually chosen to be $w_j=1$, but in many
contexts is has turned out to be useful to use weights as a kind of
cheap preconditioners. The use of \emph{adaptive weights} in sampling
theory goes back at least to~\cite{MR1247520,fhgkst}, and has become
standard in the recent work on (Fourier) sampling, 
 see for example \cite{romero,gro,gro1,strom,adnon,wfr,gataricneu,gataricneuneu,MR1882684}.

\subsection{The reconstruction operators}
In this paper we restrict ourselves to the case where the
approximation $\tilde{f} = \sum_{k\in\mathbb{N}} c_k t_k$ of the
unknown function $f\in\mathcal{H}$ is obtained by a linear and bounded 
reconstruction operator $Q:\ell
^2(\mathbb{N})\rightarrow \cT$. Thus the approximation $\tilde f$ from
the data $\langle f, u_j\rangle _{j\in \bN }$ is given by $\tilde f =
Q\big(\langle f, u_j\rangle _{j\in \bN }\big)$. 
We will  use two quantities to measure the quality of such a
reconstruction  operator. 
As a measure of stability with respect to measurement noise we use the operator norm $\|Q\|_{\rm op}$.
As a measure of stability with respect to model mismatch we follow
\cite{adhapo12} and use the so-called  quasi-optimality constant $\mu(Q)$ (see
Definition \ref{defmu}). 

Let $P_\cT$
denote the orthogonal projection onto $\cT$, $f\in \cH $ the target
function, and $l \in  \ell^2(\mathbb{N})$ be the noise vector. Then
the input data  are given by the sequence    $(\langle
f,u_j\rangle_\mathcal{H}+ l_j)_{j\in \mathbb{N}} $, the reconstruction
is $\tilde f = Q((\langle f,u_j\rangle_\mathcal{H}+ l_j)_{j\in \mathbb{N}} )$,
and the error is bounded by 
\begin{equation}\label{inequal}
\| f-Q((\langle f,u_j\rangle_\mathcal{H}+ l_j)_{j\in \mathbb{N}}
)\|_\cH\leq \mu(Q)\|f-P_\cT f\|_\cH +\|Q\|_{\rm op}\|l\|_2. 
\end{equation}

\subsection{Contributions}
The error bound \eqref{inequal} raises several questions:
\begin{itemize}
	\item Which  operators admit  an error bound of the form \eqref{inequal}?
	\item Under what circumstances does such an operator exist?
	\item Which operator has the smallest possible operator norm $\|Q\|_{\rm op}$?
	\item Which operator has the smallest possible quasi-optimality constant $\mu(Q)$?
	\item Is there a way to trade-off between quasi-optimality and operator norm?
\end{itemize} 

Our objective is to answer these questions both in 
finite-dimensional and infinite dimensional Hilbert spaces. The
results can be formulated conveniently in the language of  frame theory.

(i) \emph{Characterization of all reconstruction operators.}  We
characterize all reconstruction operators that admit an error 
estimate of the form~\eqref{inequal}. In fact, every dual frame of the set
$(P_{\cT}u_j)_{j\in \bN }$  yields a reconstruction satisfying
\eqref{inequal}. Conversely, every reconstruction operator subject to
\eqref{inequal} is the synthesis operator of a dual frame of 
$(P_{\cT}u_j)_{j\in \bN }$. Note that~\eqref{inequal} implies that 
such a reconstruction  operator $Q$ is  exact on the
reconstruction space, i.e.,  $f = Q (\langle
f,u_j\rangle_\mathcal{H})_{j\in\mathbb{N}}$ for all  $f\in \cT$. Reconstruction
operators fulfilling this property are called \textit{perfect}. 
For a precise formulation see Theorem
\ref{equivalent}.

The important insight of ~\cite{adhapo12} is the connection between 
stability  and the angle $\phi _{\cT, \cU }$ between
the sampling space and the reconstruction space. We will see that a
perfect reconstruction operator exists if and only if $\cos (\phi
_{\cT, \cU })>0$. 
It should also be mentioned that the reconstruction operators
considered in this paper are a special case of pseudoframes
\cite{pseudoframe}.

(ii) \emph{Least squares approximation.} As already mentioned, the
canonical approximation of the data $\langle f, u_j\rangle _{j\in \bN
}$ by a vector in $\cT $ is by a least squares fit. 
Let $U^*f= (\langle f,u_j\rangle_\cH)_{j\in \mathbb{N}}$ denote the
analysis operator of the frame $(u_j)_{j\in\mathbb{N}}$ and let 
$d\in \ell ^2(\mathbb{N})$ denote the vector containing the noisy
measurements $(d_j) = (\langle f,u_j\rangle_\cH  + l_j )_{j\in
  \mathbb{N}} = U^* f + l$. 
Let the reconstruction operator $Q_1$ be defined by the least squares fit 
\begin{equation}\label{lfit}
 Q_1 d = \underset{g\in \cT}{\textnormal{arg~min}} \sum_{j\in\mathbb{N}}\left|\langle g, u_j\rangle_\mathcal{H}-d_j\right|^2 
= \underset{g\in \cT}{\textnormal{arg~min}} \|U^* g - d\|_2.
\end{equation}
It is folklore that the least squares solution~\eqref{lfit} is optimal
in the absence of additional information on $f$.  Precise
formulations of this optimality were  proven in
\cite[Theorem 6.2.]{adhapo12} (including even 
non-linear reconstructions) and in 
\cite{corach} (in abstract Hilbert space). We will show in addition (Theorem
\ref{gensampl}) that   $Q_1$ is the synthesis operator
of the canonical dual frame of $(P_\cT u_j)_{j\in\mathbb{N}}$. Using
this property we derive a simple proof for the statement that the
operator $Q_1$ has the smallest possible operator norm among all
perfect reconstruction operators. 

(iii) \emph{Minimizing the quasi-optimality constant.}
Let  $W = G^{\frac{\dagger}{2}}:=(G^\dagger)^{\frac{1}{2}}$ be the
square root of the Moore-Penrose pseudoinverse of the Gramian $G =
U^*U$ of the sampling frame $\cU $ and consider  the operator $Q_0$  defined by
\begin{equation} \label{ll5}
 Q_0 d  
= \underset{g\in \cT}{\textnormal{arg~min}} \|W U^* g - W d\|_2 \, .
\end{equation}
We will show (Theorem~\ref{wichtig})  that $Q_0$ has the smallest
possible quasi-optimality constant. The reduction of the quasi-optimality
constant is one of the motivations of weighted least squares, see
\cite{gataricneu,adnon,romero,wfr}.  In \eqref{ll5} we go a step
further and use the non-diagonal matrix $W= G^{\frac{\dagger}{2}}$ as a
weight for the least squares problem. From the point of view of linear
algebra, $W$ may be seen as a preconditioner.  

In \cite{gataricneu,adnon,romero,wfr} and also
\cite{MR1882684,fhgkst,gro,gro1,gataricneuneu,MR1247520}  the
stability with respect to a bias in the measured object is considered,
i.e., the reconstruction from $U^*(f + \Delta f) = (\langle f + \Delta
f,  u_j\rangle_\mathcal{H})_{j\in \mathbb{N}}$ (stated in terms of a
frame inequality in the latter). In this context,  $Q_0$ is
the most stable operator with respect to biased objects, see the discussion in
Section \ref{stabopt}.  

(iv) \emph{Trading stability and quasi-stability.}
It is natural to ask whether one can mix between the two least squares
problems~\eqref{lfit} and \eqref{ll5}. Let
 $\Sigma_\lambda = \big(\lambda I + (1-\lambda) U^*U\big)$ and $\lambda
 \in [0,1]$ and define $Q_\lambda $ by 
\begin{equation*}
Q_\lambda d = \underset{g\in \cT}{\textnormal{arg~min}}
\|\Sigma_\lambda^{-\frac{1}{2}} U^* g - \Sigma_\lambda^{-\frac{1}{2}}
d\|_2 \, .
\end{equation*}
These reconstruction operators ``interpolate'' between $Q_1$ (most
stable with respect to noise) and $Q_0$ (most stable with respect to
model uncertainty). The parameter $\lambda $ can be seen as a
regularization parameter, or alternatively the matrix $\Sigma _\lambda
$ as version of the adaptive weights in sampling. 
In Theorem~\ref{mix} and Lemma~\ref{quasi_constant} we will study
this class of reconstruction operators and derive several
representations for $Q_\lambda $. 

(v) \emph{Fourier resampling --- numerical experiments.}
In the last part we carry out a numerical comparison of the various
reconstruction operators on the basis of the so-called resampling
problem. We approximate a function with compact support from finitely
many, nonuniform samples of its Fourier transform and then resample
the Fourier transform on a regular grid. For this problem 
we test the performance of the reconstruction operators $Q_\lambda $.

The paper is organized as follows: In Section~2  we introduce the
frame theoretic background, discuss the angle between subspaces, and
characterize all reconstruction operators satisfying the required
stability estimate \eqref{inequal}. In Section~3  we study the various
least squares problems~\eqref{lfit} and \eqref{ll5} and analyze
several representations of the corresponding reconstruction
operators. The section is complemented by general  numerical
considerations. Section~4 covers the numerical experiments on Fourier
sampling. The brief appendix collects some standard facts about
frames.

\section{Classification of all  reconstruction operators}
We will use the language of frame theory  throughout the whole paper. 
The Appendix contains a short list of basic definitions and well known
facts from frame theory. For more details on this topic, see for
instance \cite{ch08}. 

Let us introduce some notation. 
To every set of measurement  vectors $(u_j)_{j\in \mathbb{N}}$ in a Hilbert space
$\H$ (of finite or infinite dimension) we associate the synthesis
operator $U$ defined formally by $Uc= \sum _{j\in \bN} c_j u_j$ and
the \emph{sampling space}  $\mathcal{U}=\overline{\textnormal{span}}(u_j)_{j\in
  \mathbb{N}}$. The adjoint operator $U^*$  consists of the
measurements  $U^*f = ( \langle f, u_j\rangle_\mathcal{H} )_{j\in \bN }$
and is  called the analysis operator.  The frame operator is
$S = U U^*$  and the Gramian is $G = 
U^*U$. With this notation, $(u_j)_{j\in \mathbb{N}}$ is a frame for
$\mathcal{U}=\overline{\textnormal{span}}(u_j)_{j\in \mathbb{N}}$,
if there exist constants $A,B>0$, such
that for every $f\in \mathcal{U}$ 
\begin{equation*} 
  A \|f\|_\cH^2\leq \|U^* f\|_2^2\leq B\|f\|_\cH^2.
 \end{equation*}
 We always  assume that $( u_j)_{j\in \bN } $ is a
frame for $\cU  $, thus $U^*$ is bounded from $\H$ to $ \ell ^2(\mathbb{N})$ and
$U^*$ has closed range in $\ell ^2(\mathbb{N})$. We use $\R (A)$ for the range of
an operator  $A$ and $\N (A)$ for its kernel (null space). 

Likewise we assume that $(t_k)_{k\in \mathbb{N}}$ is a frame for the
\textit{reconstruction space}
$\mathcal{T}=\overline{\textnormal{span}}(t_k)_{k\in \mathbb{N}}$
with synthesis operator $T$ and analysis operator $T^*$.  Thus
\begin{equation*} 
  C \|g\|_\cH^2\leq \|T^* g\|_2^2\leq D\|g\|_\cH^2 \qquad \text{ for } g \in \cT
  \, .
\end{equation*}

Given a sequence of  linear measurements $(\langle f, u_j\rangle_\mathcal{H}
)_{j\in \bN } = U^*f$, we try to find an approximation of $f$ in the
subspace $\mathcal{T}$. Assuming that all occurring operartors are
bounded, we investigate the class of reconstruction operators $Q: \ell
^2 \to \cT $, such that $\tilde f = Q U^* f$ is the desired
reconstruction or approximation of $f$.  
We use two metrics to quantify the stability of a reconstruction  operator
${Q:\ell ^2(\mathbb{N})\rightarrow \cT}$.
As a measure for stability with respect to measurement noise we use the operator norm $\|Q\|_{\rm op}$. 
In order to measure how well $Q$ deals with the part of the function
lying outside of the reconstruction space, we use the
\textit{quasi-optimality} constant  from \cite{adhapo12}. 
\begin{defi} \label{defmu}
 Let $Q:\ell ^2(\mathbb{N})\rightarrow \cT$ and $P_\cT$ be the
 orthogonal projection onto $\cT $. 
The \textit{quasi-optimality} constant
 $\mu=\mu(Q)>0$ is the smallest number $\mu$, such that 
 \begin{equation*}
  \|f-Q U^*f\|_\cH\leq \mu \|f-P_\cT f\|_\cH,\quad \textnormal{for all } f\in \H.
 \end{equation*}
\end{defi}
If $\mu(Q) <\infty$ we call $Q$ a \textit{quasi-optimal} operator.
Since  $P_\cT f$ is the element of $\cT$ closest to $f$, the \textit{quasi-optimality} constant $\mu$
is a measure of how well $Q U^*$ performs in comparison to orthogonal
projection $P_\cT$. Note that for $f\in \cT$ we have $QU^*f = f$, thus
a quasi-optimal reconstruction operator is perfect. 

The following theorem characterizes all bounded quasi-optimal operators.

\begin{thm}\label{equivalent}
Let $\cT$ and $\mathcal{U}$ be closed subspaces of $\H$, and
$(u_j)_{j\in \mathbb{N}}$ be a Bessel sequence spanning the closed
subspace $\mathcal{U}$. For an operator $Q:\ell
^2(\mathbb{N})\rightarrow \cT$ the following are equivalent. 
\begin{enumerate}[label=(\roman*)]
\item \label{i}There exist constants $ 0\leq\mu,\beta<\infty$, such that for $f\in \H$ and $l\in \ell ^2(\mathbb{N})$ 
 \begin{equation}\label{errr1}
\| f-Q(U^*f+l)\|_\cH\leq \mu\|f-P_\cT f\|_\cH+ \beta\|l\|_2.
 \end{equation}
\item \label{ii} $Q U^*g = g$ for $g\in \cT$ and $Q$ is a bounded operator.
\item \label{iii} The sequence $(P_\cT u_j)_{j\in\mathbb{N}}$ is a frame for $\cT$.
 Let  $(h_j)_{j\in\mathbb{N}}\subset \cT$ be a dual frame of $(P_\cT
 u_j)_{j\in\mathbb{N}}$, then  $Q$ is of the form
 \begin{equation*}
  Q c = \sum_{j\in \mathbb{N}}c_j h_j,
 \end{equation*}
  i.e., $Q$ is the synthesis operator of some  dual frame of $(P_\cT u_j)_{j\in\mathbb{N}}$.
 \item \label{iiii} The operator $Q$ is bounded and $Q U^*$ is a bounded oblique projection onto $\cT$.
\end{enumerate}
\end{thm}

 Theorem~\ref{equivalent} sets up a bijection between the class of
 reconstruction operators and the class of all dual frames of $(P_\cT
 u_j )_{j\in \bN }$. 


 To prove Theorem \ref{equivalent}, we need the concept of subspace angles.
Among the  many different definitions of the angle between subspaces (see 
\cite{Szyld,st00})  the following definition is most  suitable  for our analysis.
\begin{defi}\label{defangle}
Let $\cT$ and $\mathcal{U}$ be closed subspaces of a Hilbert space
$\H$.  The subspace angle $\varphi_{\cT,\mathcal{U}}\in
[0,\frac{\pi}{2}]$ between $\cT$ and $\mathcal{U}$  is defined as  
\begin{equation}\label{angle}
\cos(\varphi_{\cT,\mathcal{U}}) = \underset{\underset{\|g\|_\cH=1}{g\in \cT}}{\inf} \|P_\mathcal{U} g\|_\cH = \underset{\underset{\|g\|_\cH=1}{g\in \cT}}{\inf}~ \underset{\underset{\|u\|_\cH=1}{u\in \mathcal{U}}}{\sup} |\langle g,u\rangle_\mathcal{H}|.
\end{equation}
\end{defi}
We observe that in general $\cos(\varphi_{\mathcal{T},\mathcal{U}})\neq
\cos(\varphi_{\mathcal{U},\mathcal{T}})$. For $\mathcal{T}\subset \mathcal{U}$,
$\cos(\varphi_{\mathcal{T},\mathcal{U}}) = 1$ and therefore
$\varphi_{\mathcal{T},\mathcal{U}} = 0$. If 
$\mathcal{U}\subsetneq \mathcal{T}$, then  $\cos(\varphi_{\mathcal{T},\mathcal{U}}) =
0$ and  $\varphi_{\mathcal{T},\mathcal{U}} = \frac{\pi}{2}$. 

The following lemma collects the main properties of oblique
projections and  angles between subspaces. 
\begin{lem}\label{gemeinsam}
Assume that  $\cT$ and $\W$ are closed subspaces of a Hilbert space
$\H$.  Then
 \begin{enumerate}[label=(\roman*)]
  \item
  \label{l}
  $\cos(\varphi_{\cT,\W^\perp})>0$ if and only if
 $\cT\cap \W=\{0\}$ and the direct sum $\cT\oplus \W$ (not necessarily orthogonal) is closed in $\H$.
 \item 
 \label{exist}
 If $\cT \cap \W = \{0\}$ and $\H_1:=\cT\oplus \W$ is a closed subspace of $\H$, then the
oblique projection $P_{\cT,\W}:\H_1\rightarrow \cT$ with range $\cT$ and kernel $\W$ is well defined and bounded on $\H_1$.
\item \label{coroll}
Let $\cos(\varphi_{\cT,\W^\perp})>0$,  $\H_1:=\cT\oplus \W$,  and let
$P_{\cT,\W}:\H_1\rightarrow \cT$ be the oblique projection with range $\cT$ and null space $\W$. Then
 \begin{equation*} 
  \|P_{\cT,\W}\|_{\rm op}=\frac{1}{\cos(\varphi_{\cT,\W^\perp})}
 \end{equation*}
 and
\begin{equation} \label{bbb}
 \|f-P_\cT f\|_\cH\leq\|f-P_{\cT,\W}f\|_\cH\leq \frac{1}{\cos(\varphi_{\cT,\W^\perp})}\|f-P_\cT f\|_\cH,
\end{equation}
for all $f\in \H_1$.
The upper bound in \textnormal{(\ref{bbb})} is sharp.
 \end{enumerate}
\end{lem}
 Item (i) of Lemma \ref{gemeinsam} is stated in \cite[Theorem
 2.1]{aha}, the proof of (ii)  can be found in \cite[Theorem
 1]{ex}, and for (iii)  see 
  \cite{Szyld}, \cite{ex}, and \cite[Corollary 3.5]{adhapo12}.

\textbf{Proof of Theorem \ref{equivalent}}
(i) $\arro $ (ii). Set $l = 0$ and choose  $f \in \mathcal{T}$. Then
\eqref{errr1} implies  $QU^* f = f$, since otherwise $\mu = \infty$.
Setting $f = 0$ in \eqref{errr1} implies  that $Q$ is bounded. 

(ii) $\arro $ (iii)  Let $Q:\ell ^2(\mathbb{N})\rightarrow \cT$ be a bounded operator with $Q U^*g = g$ for $g\in \cT$.
Let $(e_j)_{j\in \mathbb{N}}$ be the standard basis of $\ell ^2(\mathbb{N})$ and let $h_j = Q e_j$. Then 
$Q c = \sum_{j\in \mathbb{N}}c_j h_j$.
In particular for $g\in \cT$,
\begin{equation*}
Q U^* g =  \sum_{j\in \mathbb{N}}\langle g,P_\cT u_j\rangle_\cH  h_j = g.
\end{equation*}
Since $Q$ is bounded, 
$(h_j)_{j\in\mathbb{N}}$ is a Bessel sequence in $\cT$. By assumption
$(u_j)_{j\in \mathbb{N}}$ is a Bessel sequence in $\mathcal{U}$ with
Bessel bound $B$ 
and consequently
\begin{equation*}
\begin{aligned}
 \sum_{j\in \mathbb{N}} |\langle f,P_\cT u_j\rangle_\cH|^2 = \sum_{j\in \mathbb{N}} |\langle P_\cT f, u_j\rangle_\cH|^2\leq B \|P_\cT f\|_\cH^2\leq B \|f\|_\cH^2.
 \end{aligned}
\end{equation*}
Therefore $(h_j)_{j\in\mathbb{N}}$ is a dual frame of $(P_\cT u_j)_{j\in\mathbb{N}}$.


(iii) $\arro $ (iv)
Let $(h_j)_{j\in\mathbb{N}}$ be  a dual frame of $(P_\cT
u_j)_{j\in\mathbb{N}}$ and define  $Q$  by $Q c = \sum_{j\in
  \mathbb{N}}c_j h_j$ and $P:= Q U^*$. Since  the range of $P$ is
contained in $\cT $ and $Q U^* g = g$ for $g\in \cT$, it follows that
$P$ is onto $ \cT$ and that
$P^2 = Q U^* Q U^* = Q U^* = P$.
Since both $Q$ and $ U^*$ are bounded, $P$ is bounded.

(iv) $\arro $ (i).  Let $Q$ be a bounded operator, and let $P:=
Q U^*$ be a bounded oblique projection onto $\cT$. Lemma
\ref{gemeinsam}\ref{coroll}  implies that 
$\|f- P f\|_\cH \leq \|P\|_{\rm op}\|f-P_\cT f\|_\cH$,
and consequently
\begin{equation*}
 \| f-Q( U^*f+l)\|_\cH\leq \|P\|_{\rm op}\|f-P_\cT f\|_\cH + \|Q\|_{\rm op}\|l\|_2.
\end{equation*}
This finishes the proof.
\hfill $\Box$

As a direct consequence of Theorem \ref{equivalent} and  Lemma \ref{gemeinsam}, \ref{coroll}, we obtain the following characterization of the quasi-optimality constant.
\begin{cor}\label{use}
If $Q:\ell ^2(\mathbb{N})\rightarrow \cT$ is a bounded and perfect reconstruction operator, then $P = Q U^*$ is a bounded oblique projection onto $\cT$. If $\W^\perp$ denotes the null-space of $P$, then 
\begin{equation*}
 \mu(Q) = \|Q U^*\|_{\rm op} = \frac{1}{\cos(\varphi_{\cT,\W})}.
\end{equation*}
\end{cor}

In the following we always use the assumption that  the angle
between the reconstruction and sampling space fulfills
${\cos(\varphi_{\cT,\mathcal{U}})>0}$. 
The following lemma shows that this assumption is equivalent to 
$(P_\cT u_j)_{j\in \mathbb{N}}$ forming a frame for $\cT$ for every frame $(u_j)_{j\in \mathbb{N}}$ for $\mathcal{U}$.
By Theorem \ref{equivalent} \ref{iii} this is necessary for the existence of a quasi-optimal operator. In finite dimensions, for a basis $(u_j)_{j=1}^n$for $\mathcal{U}$, $(P_\cT u_j)_{j=1}^n$ can only be a spanning set for $\cT$ if $\dim(\mathcal{U}) \geq \dim(\cT)$. 
This means that by the assumption ${\cos(\varphi_{\cT,\mathcal{U}})>0}$ we restrict ourselves to an oversampled regime. 
\begin{lem}\label{proframe}
 If $\cT$ and $\mathcal{U}$ are closed subspaces of $\H$, then the following are equivalent:
 \begin{enumerate}[label=(\roman*)]
  \item $\cos(\varphi_{\cT,\mathcal{U}})>0$.\label{1111}
  \item For every frame $(u_j)_{j\in \mathbb{N}}$ for
    $\mathcal{U}$ with frame bounds $A$ and $B$, the projection
    $(P_\cT u_j)_{j\in \mathbb{N}}$ is a frame for $\cT$ with frame
    bounds $A\cos^2(\varphi_{\cT,\mathcal{U}})$ and $B$.\label{111} 

\noindent If one of these conditions is satisfied, then the following property holds:
\item $\R
(T^* U) = \R (T^*)$, therefore  both  $\R (T^*U) $ and $\R
(U^* T)$ are closed subspaces and $U^*T$ is pseudo-invertible.
Furthermore, 
\begin{equation}
  \label{eq:ll4}
  \N (U^*) \cap \cT = \{0\} \, .
\end{equation}
 \end{enumerate}
\begin{proof}
(i) $\arro $ (ii) 
Let $(u_j)_{j\in \mathbb{N}}$ be a frame for $\mathcal{U}$ with 
frame bounds $A$ and  $B$.   
 The assumption $\cos(\varphi_{\cT,\mathcal{U}})>0$ and the definition
 of $\varphi_{\cT,\mathcal{U}}$ imply  that
 \begin{equation}\label{nonu}
 \|g\|_\cH\cos(\varphi_{\cT,\mathcal{U}})\leq \|P_\mathcal{U} g\|_\cH\quad \textnormal{for all }g\in \cT.
 \end{equation}
In particular, for $g\in \cT$ we obtain with \eqref{nonu}
\begin{equation}\label{ee}
 A\|g\|_\cH^2\cos^2(\varphi_{\cT,\mathcal{U}})\leq A\|P_\mathcal{U} g\|_\cH^2
\leq \sum_{j\in \mathbb{N}}|\langle P_\mathcal{U} g,u_j\rangle_\cH |^2 \leq B \| P_\mathcal{U}g\|_\cH^2\leq B \| g\|_\cH^2. 
\end{equation}
The identity $\langle P_\mathcal{U} g,u_j\rangle_\mathcal{H}= \langle g,u_j\rangle_\mathcal{H}
= \langle g,P_\cT u_j\rangle_\mathcal{H}$ for $g\in\mathcal{T}$ and $j\in \bN $ now
shows that
$(P_\cT u_j)_{j\in \mathbb{N}}$ is a frame for $\cT$ with frame
bounds $A \cos^2(\varphi_{\cT,\mathcal{U}})$ and $B$. 

(ii) $\arro $ (i) Let $(u_j)_{j\in \mathbb{N}}$ be a frame for $\mathcal{U}$ with upper frame bound $B$ and let
$(P_\cT u_j)_{j\in \mathbb{N}}$ be a frame for $\cT$ with lower  frame bound $C_1>0$.
Since  $\langle g,P_\cT u_j\rangle_\mathcal{H} = \langle P_\mathcal{U} g,u_j\rangle_\mathcal{H}$ for $g\in \cT$, we obtain
\begin{equation*}
 C_1 \|g\|_\cH^2 \leq \sum_{j\in \mathbb{N}} |\langle g,P_\cT u_j\rangle_\cH|^2 = 
 \sum_{j\in \mathbb{N}} |\langle P_\mathcal{U} g,u_j\rangle_\cH|^2 \leq B \|P_\mathcal{U} g\|_\cH^2.
\end{equation*}
This implies that $\cos(\varphi_{\cT,\mathcal{U}})= \underset{\underset{\|g\|_\cH=1}{g\in \cT}}{\inf} \|P_\mathcal{U} g\|_\cH \geq \sqrt{\frac{C_1}{B}} >0$.

(ii) $\arro $ (iii) 
Since
$(u_j)_{j\in \mathbb{N}}$ and $(t_k)_{k\in \mathbb{N}}$ are 
Bessel sequences, both $ U^*$ and $ T$ are bounded, and
therefore  $ U^* T$ is also  bounded. 
The entries of $ U ^* T$ are given by
\begin{equation*}
 ( U^* T)(j,k) = \langle t_k,u_j\rangle_\cH =\langle t_k,P_\cT u_j\rangle_\cH,
\end{equation*}
and  $ U^* T$ is a cross-Gramian of two frames for $\cT$. 
Let $(\tilde{u}_j)_{j\in \mathbb{N}}$ be a dual frame of $(P_\cT u_j)_{j\in \mathbb{N}}$. 
Setting $c_j =\langle f,\tilde{u}_j\rangle_\cH$ we obtain, for  $f\in \cT$, 
 \begin{equation*}
   ( T^* U c)_k = \sum_{j\in \mathbb{N}}\langle
   f,\tilde{u}_j\rangle_\cH \langle P_\cT u_j,t_k\rangle_\cH = \langle
   f,t_k\rangle_\cH = (T^*f)_k. 
 \end{equation*}
It follows that 
\begin{equation}\label{zit}
\R( T^* U) = \R( T^*).
\end{equation}
Since $(t_k)_{k\in \mathbb{N}}$ is a frame for $\mathcal{T}$, 
$\R( T^*)$ is closed in $\ell ^2(\mathbb{N})$,  and so are $\R( T^* U)$ and
$\R( U^* T)$. This implies that both $T^*U$ and $U^*T$ possess a
pseudoinverse (see Appendix~\ref{pseudoinverse}). 

To prove \eqref{eq:ll4}, let $g\in \N (U^*) \cap \cT $. Then $g= Tc$ for
some $c\in \ell ^2(\bN )$ and $U^*g = U^*Tc = 0$. This means that
$c\in \N (U^*T) =   \R( T^* U)^\perp =  \R( T^*)^\perp =  \N(  T)$.
Consequently, $g =Tc = 0$, and $\N (U^*) \cap \cT = \{0\}$. 
\end{proof}
\end{lem}


\section{  The reconstruction operators}
 
\subsection{Least squares and the operator $Q_1$}\label{gens}
We first consider the reconstruction operator $Q_1:\ell ^2(\mathbb{N})\rightarrow \cT$
corresponding to the  solution of the least squares problem
\begin{equation} \label{lsrmrm}
Q_1 d = \underset{g\in \cT}{\textnormal{arg~min}}
\sum_{j\in\mathbb{N}}\left|\langle g, u_j\rangle_\cH-d_j\right|^2  
= \underset{g\in \cT}{\textnormal{arg~min}} \|U^* g - d\|_2.
\end{equation}
This approach is analyzed in detail in \cite{adhapo12}. Least square
approximation is by far the most frequent approximation method in
applications and 
of fundamental
importance, since it has the smallest operator norm among all perfect
operators. 

The following theorem reviews several representations  of
the operator $Q_1$. 
The connection of the operator $Q_1$ to the oblique projection
$P_{\cT,S(\cT)^\perp}$ was already derived in \cite[Section
4.1.]{adhapo12}  for
finite dimensional space $\mathcal{T}$. Our  new contribution 
is the connection to the canonical dual frame and the systematic
discussion of the various representions of a least squares problem. As we will apply the
statement several times,  we
include a streamlined proof. As usual, $A^\dagger $ denotes the Moore-Penrose
pseudo-inverse of an operator $A$. For the existence of $A^\dagger $
it suffices to show that the range of $A$ is closed (see Appendix~\ref{pseudoinverse}).  

\begin{thm}\label{gensampl}
Let $\cT$ and $\mathcal{U}$ be closed subspaces of a Hilbert space $\H$
such that $\cos(\varphi_{\cT,\mathcal{U}})>0$. 
Let $(u_j)_{j\in \mathbb{N}}$ be a frame for $\mathcal{U}$ with
synthesis operator $U$ and frame operator $S$. Let $(t_k)_{k\in
  \mathbb{N}}$ be a frame for $\cT$ with synthesis operator $T$. 
	
Consider the following operators:
\begin{enumerate}[label=(\roman*)]
\item $A_1=T\left(U^*T\right )^\dagger$.\label{ti}
\item The operator $A_2$ is  given  on $\R ( U ^*)$ by \label{tii}
 \begin{equation}\label{hal}
A_2  U^* = P_{\cT,S(\cT)^\perp}
 \end{equation}
 and on $\R ( U ^*)^\perp $ by
\begin{equation}\label{bb1}
 A_2c = 0 \quad\mbox{ for } c\in\R( U^*)^\perp.
\end{equation}
By \eqref{hal}  $A_2$ is independent of the particular choice of the
reconstruction frame $(t_k)_{k\in \mathbb{N}}$ for $\cT$. 
\item Let $(h_j)_{j\in \bN } $ be  the canonical dual frame of
  $(P_\cT u_j)_{j\in\mathbb{N}}$ and $A_3c = \sum _{j\in \bN } c_j
  h_j$ be the synthesis operator of $(h_j)_{j\in \bN }$. 
\item \label{tiiii}
Let  $d\in \ell ^2(\mathbb{N})$ and let 
 $\hat{c} = (\hat{c}_{k})_{k\in\mathbb{N}}$ be the unique minimal norm element of the set
\begin{equation}\label{bayern}
K:=\underset{c \in\ell^2(\mathbb{N}) }{\textnormal{arg~min}}~\| U^* T c -d\|_2.
\end{equation}
Let the operator $A_4$ be defined by $A_4 d= \sum_{k =
  1}^\infty \hat{c}_k t_k = T\hat{c} $. 
\end{enumerate}

Then all four operators are equal, $Q_1 := A_1 =A_2 = A_3 = A_4$ and provide
the unique solution to the least squares problem 	
\begin{equation*}
Q_1 d =  \underset{g\in \cT}{\textnormal{arg~min}}~ \sum_{j \in \mathbb{N}} |\langle g,u_j\rangle_\cH -d_j|^2 = \underset{g\in \cT}{\textnormal{arg~min}}
\| U^* g - d\|_2^2.
\end{equation*}

\begin{proof}
\textbf{Step 1.} First we  check that each $A_j, j= 1, \dots , 4$, is
well defined from  $ \ell
^2(\mathbb{N}) $ to $ \mathcal{T}$.  For $A_1$ this is clear by virtue of Lemma \ref{proframe}. 

 For $A_2$  we need to show that the projection $P_{\cT,S(\cT)^\perp}$
 is well defined and bounded on the whole space $\mathcal{H}$.  
According to Lemma~\ref{gemeinsam}(i)   we need to verify that $S(\cT)$ is closed,
$\cos(\varphi_{\cT,S(\cT)})>0$ and that $\H=\cT \oplus
S(\cT)^\perp$. For this we exploit the  frame inequality \eqref{ee} of $(u_j)_{j\in\mathbb{N}}$, \eqref{nonu}, and the fact
that $S = S P_\mathcal{U} =  P_\mathcal{U} S P_\mathcal{U} $, and  we obtain 
\begin{equation}\label{closes}
A \cos (\varphi_{\cT,\mathcal{U}})\|g\|_\cH \leq A \|P_\mathcal{U}
g\|_\cH \leq\|S P_\mathcal{U} g\|_\cH=\|S g\|_\cH \quad \text{for
}g\in \cT \, .
\end{equation}
The lower bound implies that  $S(\cT)$ is closed. For the angle
$\varphi_{\cT, S(\cT)}$ we obtain 
\begin{equation}\label{eds}
 \cos(\varphi_{\cT, S(\cT)})=\underset{\underset{g \neq 0}{g\in
     \cT}}{\inf}~ \underset{\underset{S h\neq 0}{h\in \cT}}{\sup}
 ~\frac{|\langle g,S h\rangle_\cH|}{\|g\|_\cH \|S h\|_\cH} \geq \underset{\underset{g \neq 0}{g\in
     \cT}}{\inf}~\frac{\langle g,S g\rangle_\cH}{\|g\|_\cH \|S
   g\|_\cH}. 
\end{equation}
Since $\langle g, Sg\rangle = \langle P_\cU g, S P_\cU g\rangle
\rangle \geq A \|P_\cU g\|_{\cH } ^2$ and $\|S g \|_{\cH } = \|S P_\cU
g \|_{\cH } \leq B \|P_\cU  g\|_{\cH }$, we continue \eqref{eds}   as
follows:  
\begin{equation}
  \label{eq:ll1}
   \cos(\varphi_{\cT, S(\cT)}) \geq  \underset{\underset{g \neq 0}{g\in
     \cT}}{\inf}~\frac{A \|P_\cU g\|_{\cH } ^2}{B \|g\|_\cH \|P_\cU  g\|_{\cH }}
 = \frac{A}{B} \,    \cos(\varphi_{\cT,\cU}) >0 \, .
\end{equation}

It remains to prove that $\cT\oplus S(\cT)^\perp = \H$, or, 
equivalently, that 
\begin{equation*}
 {(\cT\oplus S(\cT)^\perp)}^\perp = \cT^\perp\cap
 \overline{S(\cT)}=\cT^\perp \cap S(\cT) = \{0\} \, .
\end{equation*}
So assume that $g\in \cT ^\perp  \cap  S(\cT )$. Since $\R (T) = \cT$,
we may write every $t\in \cT $ as $t=Td$ for some $d\in \ell ^2(\bN
)$. In particular, there exist $c\in \ell ^2(\bN)$  and $v = Tc\in \cT
$, such that $g = Sv = STc$. Then for all $d\in \ell ^2(\bN )$,  the
element $g\in \cT ^\perp  \cap  S(\cT )$ satisfies
\begin{align*}
0= \langle g,t\rangle  = \langle STc,Td\rangle = \langle
UU^*Tc,Td\rangle = \langle U^*Tc, U^*Td\rangle \, . 
\end{align*}
Setting $d=c$, we obtain  $U^*Tc=0$. By Lemma~\ref{proframe}(iii)
$v=Tc \in \cT \cap \N (U^*) = \{0\}$, and thus $g=Sv = 0$, which implies
that $ \cT ^\perp \cap  S(\cT ) = \{0\}$.


The operator $A_3$ is the synthesis operator with respect to  the
canonical dual frame of $(P_\cT u_j)_{j\in\mathbb{N}}$ and is
therefore bounded by general frame theory. 

Now to $A_4$:  
By Lemma \ref{proframe} the operator $U^*T$  has a closed range and therefore its Moore-Penrose pseudoinverse is well defined.
It is well known that $\hat{c} = ( U^* T)^\dagger d$ is the unique
element of $K$ of minimal norm. Consequently, $A_4 d = \sum _{k\in \bN
} \hat{c}_k t_k$ is bounded on $\ell ^2(\bN )$. 


\textbf{Step 2.} We next show that all these operators are equal. 

Claim $A_1 = A_4$.  Since $\hat{c} = ( U^* T)^\dagger d$ is the unique
element of $K$ of minimal norm and $A_4d = T\hat{c} = T (U^*T)^\dagger
d$, we have $A_1 = A_4$. 

Claim $A_1 = A_2$.   
We define $R:= A_1 U^* = T( U^* T)^\dagger U^*$ and show that $R^2 = R$, $\R(R) =
\cT$ and $\N(R) = S(\cT)^\perp$. 
 The equality $R^2 = R$ follows from the identity  $A^\dagger A
 A^\dagger = A^\dagger$  for  the Moore-Penrose
 pseudoinverse applied to  $A = U^*T$. 
Clearly $\R(R)\subseteq \cT$.
To prove the converse inclusion  we show that $\R(RT) = \mathcal{T}$.
Using  $\R (T^* U ) = \R (T^*)$ from Lemma~\ref{proframe}
and  $A^\dagger A = P_{\R(A^*)}$ we conclude  that
\begin{equation*}
\begin{aligned}
 R T = T( U^* T)^\dagger  U^*T = T P_{\R(T^* U)} = T P_{\R (T^*)} = T P_{\N(T)^\perp} = T,
 \end{aligned}
\end{equation*}
which proves $\R(R) = \cT$. 

Now let $f\in \N (R)$, then  we have, for all $h\in \cH $, 
$$
\langle T(U^*T)^\dagger U^*f, h \rangle = \langle (U^*T)^\dagger
U^*f, T^* h \rangle = 0 \, .
$$
Since $\R (T^*) = \R (T^*U)$ by Lemma~\ref{proframe}(iii), this means
that for all $c\in \ell ^2(\bN )$
\begin{align*}
  0&= \langle (U^*T)^\dagger U^*f, T^* Uc \rangle  =  \langle U^*T
  (U^*T)^\dagger U^*f, c \rangle \\
& = \langle P_{\R (U^*T)} U^*f , c\rangle  = \langle f , U  P_{\R
  (U^*T)}c\rangle \\
&= \langle f , UU^*Tc \rangle  \, .
\end{align*}
In other words, $f\in \R (U^*UT)^\perp = S(\cT )^\perp$, as claimed.


Claim $A_1 = A_3$.  We need to show that the operator $A_1 = T\left(U^*T\right )^\dagger$ is the synthesis operator of the canonical dual frame of $(P_\cT u_j)_{j\in\mathbb{N}}$.
The frame operator $\tilde{S}$ of $(P_\cT u_j)_{j\in \mathbb{N}}$ can be written in the form 
\begin{equation*}
\tilde{S}f = \sum_{j\in \mathbb{N}} \langle f , P_\cT u_j\rangle_\mathcal{H} P_\cT u_j =
P_\cT\Big(\sum_{j\in \mathbb{N}} \langle P_\cT f ,  u_j\rangle_\mathcal{H}  u_j\Big)
= P_\cT  U U^*P_\cT f \, .
\end{equation*}
By Definition~\ref{frfr}(iv),    the canonical dual frame
of $(P_\cT u_j)_{j\in \mathbb{N}}$  is  given by $(\tilde S ^\dagger
P_{\cT } u_j) _{j\in \bN }$ with synthesis operator  
\begin{equation} \label{a3}
 A_3 c = \sum_{j\in \mathbb{N}} c_j (P_\cT  U U^*P_\cT)^\dagger P_\cT
 u_j = (P_\cT  U U^*P_\cT)^\dagger P_\cT U c = ( U^*P_\cT)^\dagger c, 
\end{equation}
where we used  $A^\dagger = (A^*A)^\dagger A^*$ with $A =
U^*P_\cT$ for the last equality. Since we have already proved that
$A_1=A_2$, we know that the operator $A_1$  is   independent of the
particular choice of a  frame  for $\cT$. We may therefore use the
frame  $(P_\cT u_j)_{j\in \mathbb{N}}$ with synthesis operator $P_\cT
U$ instead of $T$, and as a consequence obtain that  $A_1 = P_\cT  U ( U^* P_\cT  U)^\dagger = 
( U^*P_\cT)^\dagger$, where now we use $A^\dagger = A^*(A
A^*)^\dagger$ with $A =  U^*P_\cT$. Comparing with \eqref{a3}, we have
proved that  $A_3 = A_1$.

\textbf{Step 3.} Finally we show that each  operator $A_1 = \dots =
A_4$  provides the unique solution to the least squares fit~\eqref{lsrmrm}. 
Since $\N( U^*)\cap \cT = \{0\}$ by Lemma~\ref{proframe}(iii),
the solution $\tilde{f}\in \cT$  of the least squares problem 
\begin{equation}\label{sofaad}
 \tilde{f} = \underset{g\in \cT}{\textnormal{arg~min}}~\| U^* g- d\|_2^2
\end{equation}
is unique. 
Since $\R( T) = \cT$,  there exists a $c\in \ell ^2(\mathbb{N})$, such
that $\tilde{f} =  T c$, and 
by \eqref{sofaad} $\tilde{f} =  T c$ for every element $c\in K$
(cf. \eqref{bayern}). In particular, for  the minimal norm element
$\hat{c} = (U^*T)^\dagger d \in K$ used for the definition of the
operator $A_4$, we obtain  $\tilde f = T\hat{c} = T(U^*T)^\dagger d =
A_4d = Q_1d$. 
\end{proof}
\end{thm}

Theorem~\ref{gensampl} implies a simple proof for the statement that the operator $Q_1$ has the
smallest possible operator norm among all perfect reconstruction
operators. This has already been proven in \cite[Theorem 6.2.]{adhapo12} in
a more general setup that includes non-linear reconstruction
operators.  
\begin{thm}\label{simplet}
Let $\cT$ and $\mathcal{U}$ be two closed subspaces of a  Hilbert
space $\H$ such that  $\cos(\varphi_{\cT,\mathcal{U}})>0$. If $Q:\ell
^2(\mathbb{N})\rightarrow \cT$ is a 
perfect reconstruction operator ($QU^* g = g$ for $g\in \mathcal{T}$), then   
 \begin{equation*}
 \|Q\|_{\rm op}\geq \|Q_1\|_{\rm op}. 
 \end{equation*}
\begin{proof}
Let $Q:\ell ^2(\mathbb{N})\rightarrow \cT$ be a bounded and perfect operator.
From Theorem \ref{equivalent} we infer that $Q$ is the synthesis
operator of a dual frame of $(P_\cT u_j)_{j\in \mathbb{N}}$.  
From Lemma \ref{candual} (expansion coefficients  with respect to the
canonical dual frame have  the minimum $\ell ^2$-norm) we infer that for $g\in \cT$
\begin{equation*}
 \|Q^*g\|_2^2 = \|Q_1^*g\|_2^2 + \|Q^*  g - Q_1^*g\|_2^2\geq \|Q_1^* g\|_2^2 \,  . 
\end{equation*}
Since $Q^*$ is the analysis operator of a frame for $\mathcal{T}$,  we
have $Q^*g^\perp = Q_1^*g^\perp =0$ for
$g^\perp \in \cT^\perp$. Therefore
$\|Q^*\|_{\rm op}\geq \|{Q_1}^*\|_{\rm op}$, and consequently $\|Q\|_{\rm op}\geq \|Q_1\|_{\rm op}$. 
\end{proof}
\end{thm}

\subsection{The operator $Q_0$} \label{fram} 
In the last section we analyzed the operator $Q_1$ with the smallest
operator norm. We now introduce  and study the operator 
$Q_0$ with the smallest  quasi-optimality constant. 
In the following
we write $G^{\frac{\dagger}{2}}= (G^\dagger ) ^{1/2}$ when $G$ is a
positive operator with a pseudoinverse. \footnote{An early version of
  Theorem~\ref{wichtig} was announced in   our technical
  report https://arxiv.org/pdf/1312.1717.pdf (Theorem 2.6)} 

\begin{thm}\label{wichtig}
Let $\cT$ and $\mathcal{U}$ be closed subspaces of a Hilbert
space $\H$ such that $\cos(\varphi_{\cT,\mathcal{U}})>0$. 
Let $(u_j)_{j\in \mathbb{N}}$ be a frame for $\mathcal{U}$ with
synthesis operator $U$ and Gramian $G = U^*U$, and  let $(t_k)_{k\in
  \mathbb{N}}$ be a frame for $\cT$ with synthesis operator $T$. Consider the following operators:
\begin{enumerate}[label=(\roman*)]
\item $B_1:= T\left(G^{\frac{\dagger}{2}} U^* T\right )^\dagger
  G^{\frac{\dagger}{2}}$.\label{0i}
\item The operator $B_2$ given on $\R ( U ^*)$ by
 \begin{equation}\label{grunzf}
B_2  U^* = P_{\cT,P_\mathcal{U}(\cT)^\perp} 
 \end{equation}
 and on $\R ( U ^*)^\perp $ by 
\begin{equation}\label{grunzf1}
 {B_2}f = 0 \quad\mbox{for } f\in\R( U^*)^\perp.
\end{equation}
Consequently $B_2$ depends only on the subspace $\cT $, but not on  the
particular choice of a  frame $(t_k)_{k\in
  \mathbb{N}}$ for $\cT$. 
\item \label{0iii}
Let  $d\in \ell ^2(\mathbb{N})$ and 
   $\hat{c} = (\hat{c}_{k})_{k\in\mathbb{N}}$ be the unique minimal
 norm element of the set 
\begin{equation}\label{eqiii}
K:=\underset{c \in \ell
^2(\mathbb{N})}{\textnormal{arg~min}}~\| U^* T c -d\|_{G^\frac{\dagger}{2}}:= \underset{c \in \ell
^2(\mathbb{N})}{\textnormal{arg~min}}~\| {G^\frac{\dagger}{2}}U^* T c -{G^\frac{\dagger}{2}}d\|.
\end{equation}
Let the operator $B_3$ be defined by  $ B_3 d=  \sum_{k
  = 1}^\infty \hat{c}_k t_k$. 
\end{enumerate}
Then the operators defined by (i)-(iv) are equivalent, $Q_0 := B_1 = B_2 = B_3$
and provide the  unique
 solution of  the least squares problem 
\begin{equation*}
Q_0 d = \underset{g\in \cT}{\textnormal{arg~min}}~
\| U^* g - d\|_{G^{\frac{\dagger}{2}}}^2. 
\end{equation*}
\begin{proof}
Let $( S^\frac{\dagger}{2} u_j)_{j\in \bN }$ be the tight frame for $\cU $
associated to $( u_j)_{j \in \bN}$. 
 By  Lemma \ref{tight1} its analysis operator $L^*$ is given by 
\begin{equation}\label{ana}
 L^* = U^* S ^{\frac{\dagger}{2}}  =  G^{\frac{\dagger}{2}} U^*.
 \end{equation}

We  now apply Theorem~\ref{gensampl} to the frames $(t_k)_{k \in \mathbb{N}}$ for $\cT$ and
$(S^{\frac{\dagger}{2}} u_j)_{j\in\mathbb{N}}$ for $\mathcal{U}$. 

Since $(S^{\frac{\dagger}{2}}
u_j)_{j\in\mathbb{N}}$ is a tight frame for $\cU$, its frame operator is
$\tilde S = L L^* = P_{\cU }$. 
 As proven in Theorem~\ref{gensampl}, 
$\cH=\cT \oplus \tilde S(\cT)^\perp = \cT \oplus P_\mathcal{U}(\cT)^\perp $ and 
 the projection $P_{\cT, \tilde S (\cT )^\perp } = P_{\cT , P_\mathcal{U}(\cT)^\perp }$ is well defined and bounded.

Let us show that $B_1 = B_2$.
We set $A_1 = T (L^* T)^\dagger $.
The equivalence of the operators $A_1$ and $A_2$ of Theorem~\ref{gensampl} says that $A_1 L^* =
P_{\cT , P_\mathcal{U}(\cT)^\perp}$ on $\R (L^*)$ and $A_1 = 0$ on $\R
(L^*)^\perp $.   Consequently, with \eqref{ana}, we obtain
$$
A_1 L^* = T(L^*T)^\dagger L^* = T(G^{\frac{\dagger}{2}} U^*T)^\dagger
G^{\frac{\dagger}{2}} U^* = B_1 U^* \, ,$$
and 
$$
B_1 U^* = A_1 L^* =  P_{\cT, P_\cU
  (\cT )^\perp} \, . 
$$
In order to prove that \eqref{grunzf1} holds for $B_1$, we show that $\R( U^*)^\perp= \N\big(G^\frac{\dagger}{2}\big)$.
The set $\R( U^*)$ is closed because $(u_j)_{j\in \mathbb{N}}$ is a  frame for $\cU
$, and therefore $\R(G) = \R( U^* U) = \R( U^*)$. 
Since  $\N(A^\dagger) = \N(A^*)$, we obtain
\begin{equation*}
\begin{aligned}
 \N\big(G^\frac{\dagger}{2}\big)= \N\big(G^\dagger\big) = \N\left(G^*\right)
 =\N(G)=\R(G)^\perp = \R( U^*)^\perp.
 \end{aligned}
\end{equation*}

That the operator $B_1$ has the representation (iii) follows from the fact that the minimal norm element of $K$ (cf. \eqref{eqiii}) is obtained by the Moore-Penrose pseudoinverse, i.e. 
$\hat{c} =  (G^{\frac{\dagger}{2}} U^*T)^\dagger G^{\frac{\dagger}{2}} d$.


By Theorem~\ref{gensampl} $A_1 = T(L^*T)^\dagger $ solves the following
least squares problems: for every $\tilde d \in \ell ^2$, in
particular for $\tilde d = G^{\frac{\dagger}{2}}d$ with $d\in \ell
^2(\mathbb{N})$, the element $\tilde f = A_1 \tilde d = A_1 G^{\frac{\dagger}{2}}d
= B_1 d $ solves 
$$
\tilde f = \underset{g\in \cT}{\mathrm{arg~min}} \, \|L^* g - \tilde d\|_2 =
\underset{g\in \cT}{\mathrm{arg~min}}
\, \| G^{\frac{\dagger}{2}}U^* g -   G^{\frac{\dagger}{2}} d\|_2. 
$$
\end{proof}
\end{thm}

\begin{remark}
{\rm The approximation or reconstruction of $f$ from $U^*f$ by means of
$Q_0$ can be understood as a two-step procedure: first the input data
$U^*f$ are  preprocessed with $G^{\frac{\dagger}{2}}$,  the result is 
  $G^{\frac{\dagger}{2}}U^*f =  L^*f = (\langle f,
    S^{\frac{\dagger}{2}}u_j\rangle_\mathcal{H})_{j\in\mathbb{N}}$. Then
      $\tilde f$ is produced with $T(L^* T)^{\dagger}$, 
			which is again  the
      synthesis operator of a frame.}
\end{remark}

The next result shows that the operator $Q_0$ has the smallest
possible quasi-optimality constant.
\begin{thm}\label{mini}
Let $\cos(\varphi_{\cT,\mathcal{U}})>0$ and let $Q_0$ be defined as in Theorem \ref{wichtig}(i).
If $Q:\ell ^2(\mathbb{N})\rightarrow \cT$ is a perfect reconstruction
operator, then 
$$
\mu(Q)\geq \mu(Q_0)
$$
or, equivalently, $\|Q  U^*\|_{\rm op}\geq \|Q_0  U^*\|_{\rm op}$.
\begin{proof}
We recall that $\mu(Q)$ is the smallest $\alpha$ such that for every $f\in \H$
  \begin{equation}\label{jaja}
  \|f-Q U^*f\|_\cH \leq\alpha\|f-P_\cT f\|_\cH \, ,
 \end{equation}
and we may assume that  $\mu(Q)<\infty$.  Let $g\in \cT$ and $u^\perp
\in \cU ^\perp$. 
 Then  inequality \eqref{jaja} implies that 
 $Q U^*g=g$ and $QU^* u^\perp = 0$. 
This means that $Q U^*f=P_{\cT,\mathcal{U}^\perp}f$ for $f\in\cT\oplus \mathcal{U}^\perp$. 
 Since by Corollary \ref{use}  
 the sharp upper bound is
\begin{equation*}
\|f-Q U^*f\|_\cH\leq\frac{1}{\cos(\varphi_{\cT,\mathcal{U}})}\|f-P_\cT f\|_\cH
\qquad  \text{for } f\in \cT\oplus \mathcal{U}^\perp \, ,
\end{equation*}
we conclude that 
$\alpha \geq \frac{1}{\cos(\varphi_{\cT,\mathcal{U}})}$.

It remains to prove that $\mu(Q_0) = \frac{1}{\cos(\varphi_{\cT,\mathcal{U}})}$.
Since $Q_0 U^* = P_{\cT,P_\mathcal{U}(\cT)^\perp}$,  Corollary \ref{use}
implies that $\mu(Q_0) = \frac{1}{\cos(\varphi_{\cT, P_\mathcal{U}
    (\cT)})}$.  
We observe that 
\begin{align*}
\cos(\varphi_{\cT, P_\mathcal{U} (\cT)})&=
\underset{\underset{\|g\|_\cH=1}{g\in \cT}}{\inf}~ \underset{\underset{\|v\|_\cH=1}{v\in P_\mathcal{U} (\cT)}}{\sup} |\langle g,v\rangle_\mathcal{H}| =
\underset{\underset{\|g\|_\cH=1}{g\in \cT}}{\inf} ~\underset{\underset{\|v\|_\cH=1}{v\in P_\mathcal{U} (\cT)}}{\sup} |\langle g,P_\mathcal{U} v\rangle_\mathcal{H}|\\
&=\underset{\underset{\|g\|_\cH=1}{g\in \cT}}{\inf} ~\underset{\underset{\|v\|_\cH=1}{v\in P_\mathcal{U}(\cT)}}{\sup} |\langle P_\mathcal{U} g,v\rangle_\mathcal{H}|=
\underset{\underset{\|g\|_\cH=1}{g\in \cT}}{\inf}\|P_\mathcal{U} g\|_\cH=
\cos(\varphi_{\cT,\mathcal{U}}),
\end{align*}
using   definition \eqref{angle} for the first equality and last
equality.  Thus $\mu (Q) \geq \mu (Q_0)$. 
\end{proof}
\end{thm}

\subsection{Combinations of $Q_0$ and $Q_1$}\label{combination}
The operators $Q_0$ and $Q_1$ optimize different performance metrics,
specifically $Q_1$ is most stable with respect to noisy data, and
$Q_0$ is optimal with respect to the deviation of the target function
from the reconstruction space. It is natural to interpolate between
these two operators and to try to define mixtures $Q_\lambda $ such
that $\|Q_1\|_{\rm op}\leq\|Q_\lambda\|_{\rm op}\leq\|Q_0\|_{\rm op}$ and
 $\mu(Q_0)\leq\mu(Q_\lambda)\leq\mu(Q_1)$, $\lambda\in(0,1)$.
To do this, we procede as follows.

For $\lambda \in [0,1]$ we define
\begin{equation}\label{wlambda}
M_\lambda = \lambda I +  (1-\lambda)S_1
\end{equation}
and
\begin{equation}\label{sigmalambda}
\Sigma_\lambda = \lambda I +  (1-\lambda)G_1
\end{equation}
where $I$ denotes the identity operator on $\H$ and  on $\ell
^2(\mathbb{N})$ respectively, $S_1 = UU^*$ the frame operator and
$G_1:= U^*U$ the Gramian of the frame $(u_j)_{j\in \mathbb{N}}$ for
$\mathcal{U}$.
For $\lambda >0$  $M_\lambda $ is invertible on $\cH$, $\mathcal{U}$
is an  invariant subspace  of  $M_\lambda $ and $\Sigma
_\lambda $ is invertible on $\ell ^2(\mathbb{N})$. We now set 
\begin{equation}\label{tildeuj}
u_{\lambda,j}: = M_\lambda^{-1/2}u_j \qquad \text{ for } j\in \bN \, .
\end{equation}
The next lemma describes the properties of the new frame $(
u_{\lambda ,j})_{j\in \bN }$. 

\begin{lem}\label{arg}
Let  $(u_j)_{j\in\mathbb{N}}$ be a frame for $\mathcal{U}$ with
frame bounds $A$ and $B$. 
Fix $\lambda\in(0,1]$, and let
$(u_{\lambda,j})_{j\in \mathbb{N}}$ be defined by \eqref{tildeuj}.

Then $(u_{\lambda,j})_{j\in\mathbb{N}}$ is a frame for $\mathcal{U}$ with  frame bounds
$\frac{A}{\lambda+ (1-\lambda)A}$ and $\frac{B}{\lambda+ (1-\lambda)B}$, i.e., for every $f\in \mathcal{U}$
\begin{equation}\label{zuzeigen}
 \frac{A}{\lambda+ (1-\lambda)A}\|f\|_\cH^2
 \leq  \sum_{j\in \mathbb{N}}|\langle f,u_{\lambda,j}\rangle_\cH |^2
 \leq \frac{B}{\lambda+ (1-\lambda)B}\|f\|_\cH^2.
 \end{equation}
 Furthermore
\begin{equation}\label{zuzeigen1}
 \Sigma_\lambda^{-\frac{1}{2}} U^* =  U^*M_\lambda^{-\frac{1}{2}},
\end{equation}
i.e., the operator $\Sigma_\lambda ^{-\frac{1}{2}}  U^*$
is the analysis operator of the frame  
$(u_{\lambda,j})_{j\in\mathbb{N}}$ for $\mathcal{U}$.

\begin{proof}
Using $S_1=U U^*$ for the frame operator of $(
u_j)_{j\in \bN }$, we obtain the following:  for $f\in \mathcal{U}$
 \begin{equation}\label{zu3}
 \begin{aligned}
  \sum_{j\in \mathbb{N}}|\langle f,u_{\lambda , j}\rangle_\cH |^2  
  &= \sum_{j\in \mathbb{N}}|\langle f,M_\lambda^{-\frac{1}{2}}u_j\rangle_\cH |^2 
  = \sum_{j\in \mathbb{N}}|\langle
  M_\lambda^{-\frac{1}{2}}f,u_j\rangle_\cH |^2  \\
& = \langle M_\lambda^{-\frac{1}{2}} S_1 M_\lambda^{-\frac{1}{2}}f,f\rangle_\cH
 = \langle S_1 M_\lambda^{-1}  f,f\rangle_\cH.
  \end{aligned}
 \end{equation}
Let $f(x) = x \big( \lambda + (1-\lambda )x\big) \inv $ on $[0,\infty
)$. Then $f$ is increasing and $S_1 M_\lambda \inv = f(S_1)$.
Consider the restriction  $S_1:\mathcal{U}\rightarrow \mathcal{U}$ to the subspace $\mathcal{U}$.
 If $\sigma
(S_1) \subseteq [A,B]\subseteq (0,\infty )$, then by the spectral
theorem $\sigma (S_1 M_\lambda \inv ) = \sigma (f(S_1)) \subseteq [f(A),
f(B)]$. 
 Combining this with \eqref{zu3} implies the frame inequality \eqref{zuzeigen}.

Identity  \eqref{zuzeigen1} is proven in Lemma \ref{tight1}. 
\end{proof}
\end{lem}

\begin{thm} \label{mix}
Let $\cT$ and $\mathcal{U}$ be closed subspaces of a separable Hilbert
space $\H$ such that  $\cos(\varphi_{\cT,\mathcal{U}})>0$. Let $(u_j)_{j\in \mathbb{N}}$ be a frame for
$\mathcal{U}$ and  $(t_k)_{k\in
  \mathbb{N}}$ be a frame for $\cT$. 
For $0< \lambda \leq 1$ let  $L_\lambda$ be  the synthesis operator of the frame
$(u_{\lambda,j} = M_\lambda^{-\frac{1}{2}}u_j)_{j\in\mathbb{N}}$ for
$\mathcal{U}$  and $S_\lambda = L_\lambda L_\lambda^*$ the
corresponding frame operator. 

 Consider the following operators:
\begin{enumerate}[label=(\roman*)]
\item $C_1:=  T (\Sigma_\lambda^{-\frac{1}{2}} U^* T)^\dagger \Sigma_\lambda^{-\frac{1}{2}}$. \label{li}
\item Let the operator $C_2$ be defined on $\R ( U ^* )$  by
\begin{equation}\label{glei}
 C_2  U^* = P_{\cT,S_\lambda(\cT)^\perp}
\end{equation}
 and
\begin{equation} \label{gleis}
 C_2 f= 0 \quad \mbox{for }f \in \R( U^*)^\perp.
\end{equation}
Consequently, $C_2$ is independent of the particular choice of the reconstruction frame $(t_k)_{k\in \mathbb{N}}$ for $\cT$.
\item  \label{fadi}
For $d\in \ell ^2(\mathbb{N})$ set $C_3 d  = \sum_{k = 1}^\infty
\hat{c}_k t_k$ with
$\hat{c} =
(\hat{c}_{k})_{k\in\mathbb{N}}$ being the minimal norm element of the
set 
\begin{equation*}
K:=\underset{c\in \ell
^2(\mathbb{N})}{\textnormal{arg~min}}~\| U^* T c -d\|_\lambda:=\underset{c\in \ell
^2(\mathbb{N})}{\textnormal{arg~min}}~
\| \Sigma_\lambda^{-\frac{1}{2}} U^* T c - \Sigma_\lambda^{-\frac{1}{2}}
d\|_2^2.
\end{equation*}
\end{enumerate}

Then these operators are equal,  $Q_\lambda := C_1 = C_2 = C_3$,  and $\tilde f = Q_\lambda d$ is the
unique  solution of the least squares problem 
\begin{equation}\label{eq_lambda}
\tilde{f} =  \underset{g\in \cT}{\textnormal{arg~min}}
~\| U^* g - d\|_\lambda ^2.  
\end{equation}
\begin{proof}
We apply Theorem~\ref{gensampl} to the frames $( u_{\lambda ,j})_{j\in
    \bN } $ for $\cU $ and $(t_k)_{k\in \N}$ for $\cT $. 

Let $L_\lambda ^*$ be the analysis operator of $( u_{\lambda ,j})_{j\in
    \bN }$ and set $A_1 = T(L_\lambda ^*T)^\dagger$. Since $A_1$ has the equivalent representation (ii) of Theorem~\ref{gensampl} 
$$
A_1 L_\lambda ^* = P_{\cT , S_\lambda (\cT )^\perp } \quad \text{ on }
\R (L_\lambda ^*)
$$
and $A_1 = 0$ on $\R (L_\lambda ^*)^\perp $. This means that
$$
 P_{\cT , S_\lambda (\cT )^\perp } = A_1 L_\lambda ^* = T(L_\lambda ^* T)^\dagger L_\lambda ^* =
T(\Sigma_\lambda ^{-\frac{1}{2}} U^* T)^\dagger \Sigma_\lambda ^{-\frac{1}{2}} U^* =
C_1 U^*.
$$
Since $M _\lambda $ is invertible and $L_\lambda ^* = U^* M_\lambda ^{-1/2}$ we have $\R (L_\lambda ^*) =
\R ( U^*M _\lambda ^{-1/2}) = \R (U^*)$ and \eqref{glei} holds for $C_1$. 
To prove that $C_1 = 0$ on $\R (U^*)^\perp $, we use three
algebraic properties of kernels: If $A,B,C $ are bounded, $A$
pseudo-invertible and $C$ invertible, then $\N(A^\dagger) = \N(A^*)$,
$\N(AB) \supseteq \N (B)$ and $\N (AC) = C\inv \N (A)$. Consequently
the kernel of $C_1 $ is 
\begin{align*}
  \N (C_1 ) &= \N (T(\Sigma _\lambda ^{-\frac{1}{2}} U^* T)^\dagger \Sigma
  _\lambda ^{-\frac{1}{2}}) \\
&= \Sigma _\lambda ^{\frac{1}{2}} \N (T (\Sigma _\lambda ^{-\frac{1}{2}}  U^* T)^\dagger) 
\supseteq \Sigma _\lambda ^{\frac{1}{2}}\N  ((\Sigma _\lambda ^{-\frac{1}{2}} U^* T)^\dagger) \\
&= \Sigma _\lambda ^{\frac{1}{2}} \N  (T^*U\Sigma _\lambda ^{-\frac{1}{2}}) =
  \N (T^*U \Sigma_\lambda \inv) \\
& \supseteq  \N (U \Sigma _\lambda \inv ) =  \R (\Sigma _\lambda \inv
U^*)^\perp \\
& =\R (U^* M_\lambda \inv)^\perp  = \R (U^*)^\perp \, .
\end{align*}
Thus $C_1 = 0$ on $\R (U^*)^\perp $, which is \eqref{gleis}. 

For showing $C_1 = C_3$ and \eqref{eq_lambda} we repeat the proof of
Theorem~\ref{wichtig} verbatim.~\end{proof} 
\end{thm}

The following Lemma gives a useful upper bound on the quasi-optimality constant of the operators $Q_\lambda$.

\begin{lem} \label{quasi_constant}
Let $\cos(\varphi_{\cT,\mathcal{U}})>0$ and let $Q_\lambda$ be defined as in Theorem \ref{mix}.
Then the quasi-optimality constant $\mu(Q_\lambda) = \|Q_\lambda U^*\|_{\rm op}$ is bounded by 
\begin{equation}\label{eq47}
\|Q_\lambda U^*\|_{\rm op}\leq \frac{1}{\cos(\varphi_{\cT,\mathcal{U}})}\sqrt{\frac{B(\lambda +A(1-\lambda))}{A(\lambda + B(1-\lambda))}}.
\end{equation}
\begin{proof}
By Lemma \ref{arg}  $(u_{j,\lambda})_{j \in \mathbb{N}}$  is a frame for $\mathcal{U}$ with  frame bounds
$\frac{A}{\lambda + (1-\lambda) A}$ and $\frac{B}{\lambda +
  (1-\lambda) B}$, and synthesis operator $L_\lambda =
\Sigma_\lambda^{-\frac{1}{2}}U^*$. Using \eqref{eq:ll1} for the frame $(u_{j,\lambda})_{j \in \mathbb{N}}$ we infer \eqref{eq47}.
\end{proof}
\end{lem}
We observe that the upper bound in \eqref{eq47} is decreasing for $\lambda\rightarrow 0$, and for $\lambda = 0$ the upper bound coincides with the operator norm $\|Q_0 U^*\|_{\rm op} = \frac{1}{\cos(\varphi_{\cT,\mathcal{U}})}$.
Unfortunately we do not know yet how to obtain a meaningful bound on the operator norm of $Q_\lambda$.

\begin{remark}
{\rm In \cite{corach} the authors  consider a
regularization term (Tikhonov regularization) . The reconstruction
operators corresponding to such a regularized least squares fit do not
fulfill $Q (\langle f,u_j\rangle_\mathcal{H})_{j\in\mathbb{N}} = f$
for $f\in \cT$ (see \cite[equation (7)]{corach}), and therefore do not
belong to the class of reconstruction operators analyzed in this
paper.}  
 \end{remark}

\subsection{Numerical calculation of the coefficients}\label{calc}

We now discuss how to calculate the coefficients of the
reconstructions for finite sequences
 $(u_j)_{j=1}^n$ and $(t_k)_{k=1}^m$ in (a possibly infinite-dimensional space) $\cH$. The reconstruction vectors $(t_k)_{k=1}^m$ are assumed to be linearly independent.
Let $d\in \mathbb{C}^n$ denote the vector consisting of the noisy measurements
\begin{equation*}
 d=  [\langle f,u_1\rangle_\mathcal{H}+l_1,\dots,\langle f,u_n\rangle_\mathcal{H}+l_n]^T. 
\end{equation*}
By Theorem~\ref{mix} the approximation  $\tilde{f} = Q_\lambda d$ of $f$
is given by the linear combination
\begin{equation*}
 \tilde{f} = \sum_{k = 1}^m \hat{c}_k t_k,
\end{equation*}
with expansion coefficients 
\begin{equation}\label{lm}
\hat{c} = \underset{c\in \ell
^2(\mathbb{N})}{\textnormal{arg~min}}~\|U^*T c- d\|_{\lambda} =
\underset{c\in \ell
^2(\mathbb{N})}{\textnormal{arg~min}}~\|\Sigma_\lambda^{-\frac{1}{2}}
U^*T c- \Sigma_\lambda^{-\frac{1}{2}} d\|_2 \, . 
\end{equation}
If we formulate this least squares problem in terms of the normal equations, we have to solve
\begin{equation}
  \label{eq:lll8}
T^*U \Sigma_\lambda^{-1} U^* T\hat{c} = T^*U\Sigma_\lambda^{-1} d \, .
\end{equation}
 We observe that the cross-Gramian $ U^* T\in
\mathbb{C}^{n\times m}$ is the matrix with entries
\begin{equation*}
( U^* T)(j,k) = \langle u_j,t_k\rangle_\cH,
\end{equation*}
and the matrix $\Sigma_\lambda \in \mathbb{C}^{n\times n}$ is given by
\begin{equation*}
\Sigma_\lambda(j,k) = 
\begin{cases}
(1-\lambda)\langle u_j,u_k\rangle_\cH &\mbox{for }j\neq k,\\
\lambda + (1-\lambda)\langle u_j,u_j\rangle_\cH &\mbox{for }j= k\,.
\end{cases}
\end{equation*}



For the solution of an overdetermined least squares problem one may
use  a direct method, such  as the QR decomposition with pivoting with
 an operation count of $\mathcal{O}(n m^2)$. 
Alternatively, one may approximate the solution of \eqref{lm} up to a
given precision $\epsilon >0$ by means of iterative methods, such as
the conjugate gradient method applied to the normal equations with an 
operation count $\mathcal{O}(\log(\epsilon) n m)$.  A concrete
realization is the LSQR algorithm, see \cite{lsqr}.

The convergence of
the conjugate gradient iteration depends fundamentally on
the condition number $\kappa (R_\lambda )$ of the matrix $R_\lambda = T^*U
\Sigma_\lambda^{-1} U^* T$ in \eqref{eq:lll8} (where $\kappa (A) =
\|A\|_{\rm op}\, \|A\inv \|_{\rm op}$).  The following lemma offers an estimate for the
condition number under the additional condition that the
reconstruction space is spanned by an orthonormal set. This is a
common practice in many applications~\cite{Hrycak,1,adhapo12,ahc1,MR1882684}. 

\begin{lem}\label{cond}
 Let $\cT$ and $\mathcal{U}$ be closed subspaces of a separable Hilbert
 space $\H$ such that  $\cos(\varphi_{\cT,\mathcal{U}})>0$. Let $(u_j)_{j\in \mathbb{N}}$ be a frame for
 $\mathcal{U}$ with frame bounds $A$ and $B$, and let $(t_k)_{k\in
   \mathbb{N}}$ be an orthonormal basis for $\cT$. 
 Set $R_\lambda =  T^* U \Sigma_\lambda^{-1} U^* T$.

Then 
 \begin{equation}\label{soba}
 \kappa(R_\lambda) \leq \frac{1}{\cos^2(\varphi_{\cT,\mathcal{U}})}~\frac{B(\lambda + A(1-\lambda) )}{A(\lambda + B(1-\lambda))}.
 \end{equation}
\begin{proof}
From Lemma \ref{arg} we know that
$(u_{j,\lambda}=M_\lambda^{-\frac{1}{2}}u_j)_{j\in\mathbb{N}}$  is a
frame for $\mathcal{U}$ with  frame bounds 
$\frac{A}{\lambda+ (1-\lambda)A}$ and
$\frac{B}{\lambda+(1-\lambda)B}$, and the synthesis operator
$L_\lambda = \Sigma_\lambda^{-\frac{1}{2}} U^*$.  
Using  $\|Tc\|_\cH = \|c\|_2$ and \eqref{ee} for the frame $(u_{j,\lambda})_{j\in\mathbb{N}}$ instead of $(u_j)_{j\in\mathbb{N}}$ we infer that
\begin{equation}\label{sokoa}
\frac{A}{\lambda + (1-\lambda) A}\cos^2(\varphi_{\cT,\mathcal{U}})\| c\|_2^2\leq \|L_\lambda^* T c\|_2^2 \leq \frac{B}{\lambda + (1-\lambda) B} \|  c\|_2^2. 
\end{equation}
Since $\kappa(R_\lambda) = \kappa(L_\lambda^*
L_\lambda) =  (\|L_\lambda\|_{\rm
  op}\,\|L_\lambda^\dagger\|_{\rm op})^2 =  \kappa(L_\lambda)^2$, inequality \eqref{soba} is now a
direct consequence of \eqref{sokoa}. 
\end{proof}
\end{lem}

\begin{remark}
{\rm 
1. We observe that the bound for $\kappa (R_\lambda )$  on the right-hand side of \eqref{soba} is
increasing in $\lambda $ and 
 we expect that also 
$\kappa(R_{\lambda_1}) \leq  \kappa(R_{\lambda_2})$ for $\lambda_1
\leq \lambda_2$. This has been tested experimentally  
in Section \ref{numex}.

2. Note that for $\lambda>0$ the solution of the original  least squares problem
$\underset{c\in  \ell
^2(\mathbb{N})}{\min} \|U^*Tc - d\|_2$ and of $\underset{c\in  \ell
^2(\mathbb{N})}{\min} \|\Sigma_\lambda^{-\frac{1}{2}} U^*Tc -
\Sigma_\lambda^{-\frac{1}{2}}  d\|_2$  are distinct in general. 
 This is an important difference to classical preconditioning of square systems, where the solution of the original and the preconditioned system coincide.

3. One may interpret the introduction of $\Sigma
_\lambda ^{-1/2}$ as a form of preprocessing of the measurement vector
$d$. In most sampling problems the preprocessing is by a diagonal matrix 
 \cite{MR1247520,romero,fhgkst,gro,gro1,strom,adnon,wfr,gataricneu,gataricneuneu,MR1882684},
where the entries are called
``adaptive weights'' or ``density compensation factors''. The use of
non-diagonal matrices seems to be a new idea. 


4. The use of more general matrices for preprocessing is very
promising, but requires additional numerical considerations. To
achieve a small numerical complexity, one needs to approximate
$\Sigma_\lambda^{-1}$ by a simpler matrix $V_\lambda$ and then  solve the normal equations
\begin{equation*}
T^* U V_\lambda U^*T c = T^* U V_\lambda d.
\end{equation*}
This question will be pursued in future work. }
\end{remark}



\subsection{Conditions for the approximations to coincide}\label{concoinside}
While in general the reconstruction operators
$Q_1$, $Q_0$ and $Q_\lambda$ are  different, they coincide in several situations. 
\begin{lem}\label{kkurz}
 Let $\cT$ and $\mathcal{U}$ be closed subspaces of a separable Hilbert
 space $\H$ and let $(u_j)_{j=1}^n$ be a frame for $\mathcal{U}$. If
 $\cT\oplus \mathcal{U}^\perp = \H$ and $Q:\ell
 ^2(\mathbb{N})\rightarrow \mathcal{T}$  is a  bounded,  perfect
 reconstruction  operator, then  
\begin{equation}\label{sowass}
 Q  U^* = P_{\cT,\mathcal{U}^\perp} \, .
\end{equation}
 Consequently for $\lambda\in[0,1]$
\begin{equation} \label{hji}
 Q_0 = Q_1 = Q_\lambda.
\end{equation}
\begin{proof}
 As in the proof of Theorem \ref{mini} we see  that $Q U^* g = g $ for
 $g\in \cT$  and $Q U^* u^\perp = 0$ for $u^\perp \in \cU ^\perp$
 imply  that $\R (QU^*) \supseteq \cT $ and $\N (QU^*) \supseteq \cU
 ^\perp $.  Since by assumption $\cT\oplus \mathcal{U}^\perp = \H$, this proves \eqref{sowass}.


Since  $Q_1 c = Q_0 c = Q_\lambda c = 0$ for $c\in\R( U^*)^\perp$,
this implies \eqref{hji}.
%
	%
\end{proof}
\end{lem}

The decomposition $\cT\oplus\mathcal{U}^\perp = \H$ is the general assumption for consistent sampling \cite{23,24,28,oleyonina,oleole1,oleole}.
In finite dimensions the assumption  $\cT\oplus
\mathcal{U}^\perp = \H$ is fulfilled only  if $\dim(\mathcal{T}) =
\dim(\mathcal{U})$ and $\cos(\varphi_{\cT,\mathcal{U}})>0$ \cite[Lemma 3.7]{adhapo12}. In case of
linearly independent sampling and reconstruction vectors, the condition
$\dim(\mathcal{T}) = \dim(\mathcal{U})$ requires as many sampling as
reconstruction vectors. 
In other words, 
in
the critical case (between overdetermined and underdetermined) all
reconstruction operators coincide.

\begin{thm}\label{rizzz1}
Let $\cT$ and $\mathcal{U}$ be closed subspaces of $\H$ such that  $\cos(\varphi_{\cT,\mathcal{U}})>0$. If $(u_j)_{j\in \mathbb{N}}$ is a tight frame for $\mathcal{U}$, then for $\lambda \in [0,1]$
\begin{equation*}
 Q_0 = Q_1 = Q_\lambda.
\end{equation*}
\begin{proof}
Since $Q_0 = Q_1 = Q_\lambda = 0$ on $\R( U^*)^\perp$, it is sufficient to show that $Q_1  U^* = Q_0  U^* = Q_\lambda U^*$.
Using the frame operator $S_\lambda$ of the frame  
$(u_{\lambda,j} = M_\lambda^{-\frac{1}{2}}u_j )_{j\in\mathbb{N}}$
 of  $\mathcal{U}$ (cf. \eqref{wlambda}), we have $Q_\lambda  U^* =
P_{\cT, S_\lambda(\cT)^\perp}$ (by
  Theorem \ref{mix}).

Since $(u_j)_{j\in \bN }$ is a tight frame for $\cU $, its frame
operator is $S_1= A P_\cU$ for some $A>0$. Consequently,  
$$M_\lambda P_\cU  = (\lambda I +  (1-\lambda)S_1) P_\cU  = (\lambda  +
(1-\lambda) A) P_\cU \, ,$$
and $M_\lambda ^{-1/2} u_j = (\lambda  +
(1-\lambda) A)^{-1/2} u_j$ is just a constant multiple of the original
tight frame.  Therefore $(M_\lambda ^{-1/2} u_j)_{j\in \bN }$ is again
a tight frame for every $\lambda \in [0,1]$, $S_\lambda (\cT
) = P_\cU (\cT)$ and $Q_\lambda U^* = P_{\cT,S_\lambda(\cT)^\perp} = P_{\cT, P_\cU (\cT )^\perp }$  is independent of $\lambda $. Consequently, $Q_0 =
  Q_\lambda = Q_1$. 

\end{proof}
\end{thm}

\subsection{Stability with respect to a biased objects}\label{stabopt}
In \cite{gataricneu,adnon,romero,wfr} and also
\cite{fhgkst,gro,gro1,gataricneuneu,MR1247520,MR1882684}  a notion of
stability with respect to 
a bias in the measured object is considered (in the latter stated in terms of a frame inequality). This means that the measurements are
made on the vector $f+\Delta f$ instead of the correct  $f$, and
$\Delta f\in \cH $ is the bias or object uncertainty.  In this case 
the error estimate is of the form
\begin{equation}\label{gerrneu}
 \|f-Q U^*(f+\Delta f)\|_\cH\leq \mu(Q)\|f-P_\cT f\|_\cH + \|Q U^*\|_{
\rm op} \, \|\Delta f\|_\cH
 \, .
\end{equation}
It is  important to understand the conceptual  difference between
\eqref{gerrneu} and \eqref{errr1}. The error estimate \eqref{errr1} treats the
error arising  from perturbed or noisy  measurements $ U^*f+l$. Estimate
\eqref{gerrneu} treats the uncertainty of the target function
(\emph{object uncertainty}) and assumes that the exact measurements of the biased function 
$f+\Delta f$ are available. 
Since the operator $Q_0$ has the smallest possible quasi-optimality
constant $\mu (Q_0)$ and since $\mu (Q_0) = \|Q_0 U^* \|_{\rm op}$,
Theorem~\ref{mini} yields  the following corollary.
\begin{cor}\label{tbvc}
 Let $\cos(\varphi_{\cT,\mathcal{U}})>0$ and let $Q_0$ be defined as in Theorem \ref{wichtig}. If an operator $Q:\ell ^2(\mathbb{N})\rightarrow \cT$ satisfies for $f \in \H$ and $\Delta f \in \H$ 
\begin{equation} \label{bay}
 \|h-Q U^*(f+\Delta f)\|_\cH\leq \beta_1 \|f-P_\cT f\|_\cH + \beta_2\|\Delta f\|_\cH
\end{equation}
for some $0<\beta_i<\infty$, then $\beta_i \geq\mu(Q_0)$, $i =1,2$.
\end{cor}
Consequently, if we restrict ourselves to linear
mappings, Corollary \ref{tbvc} shows that $Q_0$ is optimal for the problem considered in  
\cite{gataricneu,adnon,romero,wfr,fhgkst,gro,gro1,gataricneuneu,MR1247520,MR1882684}

\section{Numerical experiments for reconstruction from Fourier measurements} \label{numex}
In this section, we apply  the various reconstruction
methods to the 
reconstruction of  a compactly supported function from non-uniform
Fourier  samples. This approximation problem occurs in numerous
applications, for  example, radial sampling of the Fourier transform
is used in MRI and CT, see \cite{lewitt}.

From the given data $\hat{f}(\omega _j), j =-n , \dots , n,$ of a
compactly supported function, we
calculate the  Fourier coefficients $\hat{f}(k), k=-m , \dots , m$,  of
$f$ and construct a final approximation by a truncated Fourier
series. This is  the uniform resampling problem, see
\cite{gelbgelb,adhapo12}. If $f$ is  smooth and periodic,  then the   
Fourier series converges exponentially fast. However, if $f$ is non-periodic or  discontinuous, then the  Fourier series of $f$ converges
slowly and also  suffers from the Gibbs phenomenon. 
Of course, for discontinuous or non-periodic functions the trigonometric polynomials
of fixed degree are a bad choice for the reconstruction space. 
Since the function $f$ is unknown there will always be some model mismatch in practice, independent of the particular choice of the reconstruction vectors. Our objective in this section is not to choose optimal reconstruction functions, but rather 
to compare how the various  reconstruction
operators $Q_\lambda$ deal with the model mismatch in noisy regimes. We will see that a smart choice of the
parameter $\lambda $ yields better approximations than the standard
least square approximation~\eqref{lfit}.

We remark that the Gibbs phenomenon can be avoided  by choosing a more
appropriate reconstruction space, e.g., algebraic polynomials~\cite{Hrycak,adcock}
or wavelet expansions~\cite{1,ahc1}. 

\subsection{Setup}
We denote by $\langle f,g \rangle_{L^2} =
\int_{-\infty}^\infty \! f(x)\overline{g(x)} \, \mathrm{d}x$ the
standard inner product on the Hilbert space $L^2(\mathbb{R})$  
and 
the Fourier transform $\mathcal{F} $  on $L^2(\mathbb{R})$ with normalization  
\begin{equation*}
\mathcal{F}f(\xi) = \int_{-\infty}^\infty \! f(x)e^{- 2\pi i x \xi} \, \mathrm{d}x.
\end{equation*}
Let $\H$ be the subspace of $L^2(\mathbb{R})$ of functions with support in the interval $[-\frac{1}{2},\frac{1}{2}]$, i.e.,
\begin{equation*}
 \H = \left\{ f\in L^2(\mathbb{R}):\textnormal{supp}(f)\subset \Big[-\frac{1}{2},\frac{1}{2}\Big]\right\}.
\end{equation*}
The given data are finitely many (non-uniform) noisy Fourier measurements  
\begin{equation} \label{noisy_samples}
 d_j =\mathcal{F}f(\omega_j) + l_j,\quad j =
 -n,\dots,n \, . 
\end{equation}
where $l_j \in \mathbb{C}$ is additive noise. The  noise $l_j$ is
assumed to be i.i.d.\ Gaussian with variance 
(average power) $\sigma_{\ell }^2 $ and \ac{SNR}  $\mathrm{SNR} =
  \frac{\|f\|_{L^2}^2}{\sigma_\ell ^2(2n+1)}$.

The sampling space consists of the exponential functions
\begin{equation*}
u_j(x) = \e^{2 \pi i \omega_j x}
\chi_{[-\frac{1}{2},\frac{1}{2}]}(x),\quad j = -n,\dots,n \, ,
\end{equation*}
so that indeed $d_j =\mathcal{F}f(\omega_j) + l_j = \langle f,
u_j\rangle_{L^2} + l_j$.

The sampling frequencies  $\omega_j\in\mathbb{R}$ are chosen
\begin{equation}\label{omegaj}
 \omega_j = \frac{j}{2} +\delta_j, \quad j = -n,\dots,n,
\end{equation} 
with $\delta_j\in [-2,2]$ i.i.d.\ and  uniformly distributed over  the interval $[-2,2]$. 
The  reconstruction space for the resampling problem is spanned by 
the complex exponentials 
\begin{equation*}
t_k(x) = \e^{2 \pi i k x}\chi_{[-\frac{1}{2},\frac{1}{2}]}(x), \quad k = -m,\dots,m
\end{equation*}
with $m\leq n$. In the numerical simulations  we approximate  the 
 exponential function 
\begin{equation*}
f(x) = \e^x\chi_{[-\frac{1}{2},\frac{1}{2}]}(x),
\end{equation*}
from the noisy  Fourier measurements \eqref{noisy_samples}.
For the reconstruction we use the operators $Q_\lambda$ of  Theorem
\ref{mix}. This means that the vector $\hat{c} =
[\hat{c}_{-m},\dots,\hat{c}_m]^T$ containing the coefficients of the
approximation  
\begin{equation*}
\tilde{f}=\sum_{k = -m}^m \hat{c}_k t_k(x)
\end{equation*}
of $f$ is the solution of the least squares problem 
\begin{equation}\label{leastsquares}
 \hat{c} = \underset{c\in \ell
^2(\mathbb{N})}{\textnormal{arg~min}}~\|\Sigma_\lambda^{-\frac{1}{2}} U^*T c - \Sigma_\lambda^{-\frac{1}{2}}d\|_2^2.
\end{equation}
For the particular bases $(u_j), (t_k)$ consisting of
exponentials,   the cross-Gramian $U^*T\in \mathbb{C}^{(2n+1)\times 
  (2m+1)}$ has the entries 
\begin{equation*}
  ( U^* T)(j,k) = \langle u_j,t_k\rangle_{L^2} = \frac{\sin(\pi(\omega_j
    -k))}{\pi(\omega_j -k)} = \sinc(\omega_j -k), 
\end{equation*}
and the preconditioning matrix $\Sigma _\lambda $ is given by the
entries 
\begin{equation*}
 \Sigma_\lambda(j,k) = \big(\lambda I_{2n+1} +  (1-\lambda)G\big)(j,k) =
\begin{cases}
(1-\lambda)\sinc(\omega_j -\omega_k) &\mbox{for }j\neq k,\\
1  &\mbox{for }j = k \, .
\end{cases}
\end{equation*}

All results in this section have been averaged over $1000$ independent realizations of the sampling frequencies and the noise.

\subsection{Noisy samples}
In the first experiment we study the influence of the sampling rate
$\frac{2m +1}{2n+1}$ and the \ac{SNR} on the recovery performance of
the operators $Q_\lambda$. We approximate the exponential function 
$f(x) = \e^x\chi_{[-\frac{1}{2},\frac{1}{2}]}(x)$ from $181$  noisy
Fourier samples ($n=90$) and reconstruct in a space of trigonometric
polynomials of degree $m=10,20,30,40$ (with dimension $2m+1$).

Table \ref{table1} lists the operator norm $\|Q_\lambda\|_{\rm op}$, the
quasi-optimality constant $\mu(Q_\lambda)$, the angle
$\varphi_{\mathcal{T},\mathcal{U}}$, the condition number
$\kappa(\Sigma_\lambda^{-\frac{1}{2}} U^* T)$ of the matrix of the
least squares problem \eqref{leastsquares} 
and  the relative error $\frac{\|Q_\lambda
  h-f\|_{L^2}}{\|f\|_{L^2}}$ for $\text{SNR} = \infty$, $\text{SNR} = 20$dB and
$\text{SNR} = 10$dB 
and $m= 10$ in (a), $m = 20$ in (b) $m= 30$ in (c) and $m = 40$ in (d).
All values are listed in the form $E \pm \sigma$ where $E$ is the expected value and $\sigma$ the standard deviation of the quantity, and rounded  to the third decimal place. 
\begin{table}[h!]
\centering
\Small
\begin{subtable}{\linewidth}\centering
\tiny
\begin{tabular}{|c|  c c  c  c c c|}
\toprule
 &\multicolumn{3}{c|}{relative error $\frac{\|Q_\lambda h-f\|_{L^2}}{\|f\|_{L^2}}$}  &  & &\\ 
 & $\text{SNR}\! = \!\infty$ & $\text{SNR}\! = \! 20$dB    &\multicolumn{1}{c|} {$\text{SNR}\! = \! 10$dB}  &$\|Q_\lambda\|_{\rm op}$ &$\mu(Q_\lambda)$ &$\kappa(\Sigma_\lambda^{-\frac{1}{2}} U^*T)$\\
\midrule
		$Q_{0}$	  &$0.067 \pm    0.000$    &$0.088   \pm   0.043$      &$0.184    \pm    0.165$      &$8.288   \pm   10.884$           &$1.000   \pm 0.000$         &$1.000  \pm  0.000$\\
		$Q_{0.1}$	&$0.068  \pm    0.005$    &$0.075  \pm   0.011$      &$0.118   \pm    0.045$      &$4.024   \pm   3.907$      &$1.078   \pm    0.221$    &$1.747    \pm    1.233$\\
		$Q_{0.2}$	&$0.069  \pm    0.010$    &$0.076   \pm   0.013$      &$0.118    \pm   0.045$      &$3.936   \pm   3.807$      &$1.160   \pm    0.382$    &$2.256   \pm     1.842$\\
		$Q_{0.3}$	&$0.071   \pm    0.014$    &$0.077   \pm   0.016$      &$0.118   \pm    0.046$      &$3.880  \pm    3.744$      &$1.245    \pm   0.525$    &$2.724     \pm    2.369$\\
		$Q_{0.4}$	&$0.073   \pm     0.019$    &$0.079    \pm   0.020$      &$0.119   \pm   0.048$      &$3.838   \pm   3.699$      &$1.336    \pm   0.662$    &$3.190    \pm    2.881$\\
		$Q_{0.5}$	&$0.076  \pm     0.024$     &$0.081   \pm   0.024$      &$0.120     \pm  0.050$      &$3.804    \pm  3.662$       &$1.434   \pm    0.802$     &$3.682    \pm    3.412$\\
		$Q_{0.6}$	&$0.079   \pm    0.028$    &$0.084    \pm  0.028$      &$0.121    \pm    0.053$      &$3.776   \pm   3.632$      &$1.544    \pm   0.950$    &$4.224   \pm     3.996$\\
		$Q_{0.7}$	&$0.082   \pm    0.033$    &$0.086   \pm   0.033$      &$0.123   \pm    0.056$      &$3.752   \pm   3.608$      &$1.670    \pm    1.112$    &$4.848   \pm     4.672$\\
		$Q_{0.8}$	&$0.085   \pm    0.039$    &$0.089   \pm   0.038$      &$0.125    \pm   0.060$      &$3.733    \pm  3.587$      &$1.820   \pm     1.298$    &$5.603    \pm    5.503$\\
		$Q_{0.9}$	&$0.089  \pm     0.045$    &$0.093   \pm   0.044$      &$0.127   \pm    0.065$      &$3.718   \pm   3.572$      &$2.011    \pm    1.523$    &$6.581   \pm     6.616$\\
		$Q_{1}$	  &$0.094    \pm   0.052$      &$0.098   \pm   0.051$       &$0.131      \pm 0.071$      &$3.712   \pm   3.566$      &$2.276    \pm    1.824$    &$7.982     \pm    8.319$\\
\bottomrule
\end{tabular}
\caption{$m = 10$, $\varphi_{\mathcal{T},\mathcal{U}} = 1.7734e-08 \pm  1.3455e-08$}
\end{subtable}
\begin{subtable}{\linewidth}\centering
\tiny
\begin{tabular}{|c | c c c c c c|}
\toprule
 &\multicolumn{3}{c|}{relative error $\frac{\|Q_\lambda h-f\|_{L^2}}{\|f\|_{L^2}}$}  &  & &\\ 
 & $\text{SNR}\! = \!\infty$ & $\text{SNR}\! = \! 20$dB    &\multicolumn{1}{c|} {$\text{SNR}\! = \! 10$dB}  &$\|Q_\lambda\|_{\rm op}$ &$\mu(Q_\lambda)$ &$\kappa(\Sigma_\lambda^{-\frac{1}{2}} U^*T)$\\
\midrule
		$Q_{0}$	  &$0.048 \pm    0.000$    &$0.096  \pm    0.060$      &$0.259    \pm    0.201$      &$11.927  \pm    12.724$   &$1.000   \pm 0.000$         &$1.000  \pm  0.000$\\
		$Q_{0.1}$	&$0.049    \pm  0.003$    &$0.066    \pm   0.017$      &$0.147    \pm   0.056$      &$5.330   \pm   4.392$      &$1.097   \pm    0.201$    &$2.103    \pm    1.391$\\
		$Q_{0.2}$	&$0.050   \pm   0.006$    &$0.067   \pm   0.018$      &$0.146    \pm   0.055$      &$5.249   \pm   4.320$      &$1.194   \pm    0.349$    &$2.833    \pm    2.088$\\
		$Q_{0.3}$	&$0.052   \pm   0.009$    &$0.068  \pm    0.019$      &$0.146    \pm   0.055$       &$5.199   \pm   4.277$      &$1.293   \pm    0.480$    &$3.499   \pm     2.694$\\
		$Q_{0.4}$	&$0.054   \pm    0.012$    &$0.069  \pm    0.021$       &$0.146   \pm    0.055$      &$5.162   \pm   4.245$      &$1.398   \pm    0.607$    &$4.164   \pm     3.286$\\
		$Q_{0.5}$	&$0.056   \pm    0.015$    &$0.070  \pm    0.022$      &$0.146   \pm    0.056$      &$5.133    \pm  4.221$      &$1.510    \pm   0.735$    &$4.867   \pm     3.906$\\
		$Q_{0.6}$	&$0.058  \pm    0.018$    &$0.072   \pm   0.024$      &$0.147   \pm    0.057$      &$5.110   \pm   4.201$      &$1.634    \pm   0.869$    &$5.644   \pm     4.591$\\
		$Q_{0.7}$	&$0.061    \pm   0.021$    &$0.074   \pm   0.026$      &$0.148   \pm    0.057$      &$5.090  \pm    4.184$       &$1.776   \pm     1.014$    &$6.546   \pm     5.387$\\
		$Q_{0.8}$	&$0.064   \pm    0.025$    &$0.076  \pm    0.029$      &$0.148   \pm    0.059$      &$5.074  \pm    4.171$      &$1.945    \pm    1.180$    &$7.650   \pm     6.372$\\
		$Q_{0.9}$	&$0.067   \pm    0.028$    &$0.079   \pm    0.031$      &$0.150   \pm    0.060$      &$5.062   \pm   4.162$      &$2.157    \pm    1.378$     &$9.102    \pm    7.695$\\
		$Q_{1}$	  &$0.071   \pm    0.032$    &$0.082   \pm   0.035$      &$0.151   \pm    0.063$      &$5.057   \pm   4.158$      &$2.451    \pm    1.638$    &$11.233    \pm    9.726$\\
\bottomrule
\end{tabular}
\caption{$m = 20$, $\varphi_{\mathcal{T},\mathcal{U}} = 2.0805e-08 \pm 1.2298e-08$}
\end{subtable}
\begin{subtable}{\linewidth}\centering
\tiny
\begin{tabular}{|c | c c  c  c c c|}
\toprule
 &\multicolumn{3}{c|}{relative error $\frac{\|Q_\lambda h-f\|_{L^2}}{\|f\|_{L^2}}$}  &  & &\\ 
 & $\text{SNR}\! = \!\infty$ & $\text{SNR}\! = \! 20$dB    &\multicolumn{1}{c|} {$\text{SNR}\! = \! 10$dB}  &$\|Q_\lambda\|_{\rm op}$ &$\mu(Q_\lambda)$ &$\kappa(\Sigma_\lambda^{-\frac{1}{2}} U^*T)$\\
\midrule
		$Q_{0}$	  &$0.039  \pm   0.000$     &$0.122    \pm   0.084$      &$0.373    \pm    0.292$      &$17.371  \pm    17.536$   & $1.000  \pm  0.000 $   &$1.000 \pm   0.000$\\
		$Q_{0.1}$	&$0.040   \pm   0.003$    &$0.069   \pm   0.021$      &$0.180    \pm   0.066$      &$6.452  \pm    4.861$     & $1.129    \pm   0.213 $  & $2.431     \pm   1.552$\\
		$Q_{0.2}$	&$0.042   \pm   0.006$     &$0.069   \pm   0.021$      &$0.180   \pm    0.066$      &$6.369   \pm   4.793$    &  $1.246   \pm    0.373 $  & $3.348     \pm   2.324$\\
		$Q_{0.3}$	& $0.044   \pm   0.008$    &$0.070   \pm    0.022$      &$0.178  \pm     0.066$     & $6.318  \pm    4.752$    &  $1.363   \pm    0.514 $  & $4.180     \pm   2.998$\\
		$Q_{0.4}$	&$0.046   \pm   0.011$    &$0.071     \pm  0.023$      &$0.178    \pm   0.066$     &  $6.282   \pm   4.723$    &  $1.482    \pm   0.649 $  &$ 5.008    \pm    3.658$\\
		$Q_{0.5}$	&$0.048    \pm   0.013$    &$0.072   \pm   0.023$     & $0.178   \pm    0.066$     &  $6.254  \pm    4.701$    &   $1.610   \pm     0.785 $  &$ 5.885     \pm   4.350$\\
		$Q_{0.6}$	&$0.050    \pm   0.015$    &$0.073   \pm   0.024$     & $0.179   \pm    0.066$     & $6.232  \pm    4.683$    &  $1.745    \pm   0.928 $ & $ 6.856     \pm   5.115$\\
		$Q_{0.7}$	&$0.052    \pm   0.017$    &$0.074   \pm   0.025$     & $0.179   \pm    0.067$     & $6.213   \pm   4.667$    &  $1.900   \pm     1.085 $ & $ 7.987    \pm    6.007$\\
		$Q_{0.8}$	&$0.054    \pm   0.020$    &$0.075    \pm    0.027$     & $0.179   \pm    0.068$     & $6.198  \pm    4.655$    &  $2.082   \pm     1.265 $  &$ 9.378    \pm    7.113$\\
		$Q_{0.9}$	&$0.056    \pm    0.022$    &$0.077   \pm   0.028$     & $0.180    \pm   0.068$     & $6.187   \pm   4.646$    &  $2.310   \pm     1.483 $ & $ 11.223    \pm       8.600$\\
		$Q_{1}$	  &$0.059    \pm   0.024$    &$0.079   \pm   0.030$      & $0.181    \pm    0.069$     & $6.183   \pm   4.642$    &  $2.621   \pm      1.775 $ & $ 13.966    \pm    10.885$\\
\bottomrule
\end{tabular}
\caption{$m = 30$, $\varphi_{\mathcal{T},\mathcal{U}} = 3.5152e-08 \pm  7.3768e-09$}
\end{subtable}
\begin{subtable}{\linewidth}\centering
\tiny
\begin{tabular}{|c | c c  c  c c c|}
\toprule
 &\multicolumn{3}{c|}{relative error $\frac{\|Q_\lambda h-f\|_{L^2}}{\|f\|_{L^2}}$}  &  & &\\ 
 & $\text{SNR}\! = \!\infty$ & $\text{SNR}\! = \! 20$dB    &\multicolumn{1}{c|} {$\text{SNR}\! = \! 10$dB}  &$\|Q_\lambda\|$ &$\mu(Q_\lambda)$ &$\kappa(\Sigma_\lambda^{-\frac{1}{2}} U^*T)$\\
\midrule
		$Q_{0}$	  &$0.034   \pm  0.000$&$   0.870     \pm   3.781    $&$     2.550    \pm     9.917   $&$   180.730   \pm   790.710    $&$       1.000  \pm  1.000    $&$     1.000   \pm 1.000$\\
		$Q_{0.1}$	&$0.036   \pm   0.004$&$   0.075   \pm   0.025   $&$   0.209    \pm   0.076  $&$    7.322   \pm   5.272   $&$    1.185    \pm   0.298   $&$  2.701   \pm     1.705$\\
		$Q_{0.2}$	&$0.038   \pm   0.007$&$   0.075   \pm   0.025   $&$   0.208    \pm   0.073   $&$   7.229   \pm   5.194   $&$   1.327    \pm   0.469  $&$  3.758    \pm    2.537$\\
		$Q_{0.3}$	&$0.040   \pm   0.009$&$    0.076   \pm   0.026   $&$   0.207    \pm   0.072   $&$   7.179    \pm  5.153   $&$   1.456   \pm    0.610  $&$  4.713    \pm    3.263$\\
		$Q_{0.4}$	&$0.041   \pm    0.011$&$   0.076   \pm   0.026   $&$   0.206    \pm   0.072   $&$   7.146   \pm   5.127   $&$   1.580   \pm    0.739  $&$  5.667    \pm    3.978$\\
		$Q_{0.5}$	&$0.043    \pm    0.012$&$   0.077   \pm   0.027   $&$   0.206    \pm   0.071   $&$   7.121   \pm   5.107   $&$   1.707   \pm    0.865  $&$  6.677   \pm     4.729$\\
		$Q_{0.6}$	 &$0.044   \pm    0.014$&$   0.077  \pm    0.027  $&$    0.206    \pm   0.071   $&$   7.102   \pm   5.093   $&$   1.840   \pm    0.993  $&$  7.800    \pm    5.562$\\
		$Q_{0.7}$	&$0.046    \pm    0.015$&$    0.078   \pm   0.028    $&$  0.206     \pm   0.071   $&$   7.088    \pm  5.081    $&$  1.983    \pm    1.129   $&$   9.110    \pm    6.535$\\
		$Q_{0.8}$	&$0.048   \pm    0.017$&$    0.079  \pm    0.028   $&$   0.206   \pm    0.071   $&$   7.076   \pm   5.073   $&$   2.147   \pm     1.281  $&$  10.729    \pm    7.744$\\
		$Q_{0.9}$	&$0.049   \pm    0.018 $&$    0.080  \pm    0.029   $&$   0.206    \pm   0.071   $&$   7.068   \pm   5.067   $&$   2.344     \pm     1.460  $&$  12.888    \pm    9.374$\\
		$Q_{1}$	  &$0.051    \pm   0.019$&$     0.081   \pm   0.029   $&$   0.207    \pm   0.072   $&$   7.065   \pm   5.065   $&$   2.601    \pm    1.689  $& $  16.127    \pm    11.882$\\
\bottomrule
\end{tabular}
\caption{$m = 40$, $\varphi_{\mathcal{T},\mathcal{U}} = 5.4254e-08 \pm   1.4382e-08$}
\end{subtable}
		\caption{Reconstruction of the exponential function by
                  the 
                  operators $Q_\lambda$ from noisy
                  measurements}\label{table1} 
\end{table}

By  \eqref{inequal} the (absolute) reconstruction error depends
both  on the quasi-optimality constant $\mu(Q_\lambda)$ and the operator
norm $\|Q_\lambda\|_{\rm op}$. 
The numerical simulations support Theorem \ref{mini} asserting that
the  quasi-optimality constant
$\mu$ is minimal for $\lambda =0$. As expected in view of Lemma \ref{quasi_constant} the quasi-optimality constant
$\mu$ is increasing with $\lambda$.  
The angle $\varphi_{\mathcal{T},\mathcal{U}}$ is almost zero, so the reconstruction space is
  ``almost contained'' in the sampling space. Therefore $Q_0 U^*$ is nearly identical to the orthogonal projection $P_\mathcal{T}$ onto $\mathcal{T}$.
A small angle $\varphi_{\mathcal{T},\mathcal{U}}$ is essential for stable reconstruction. Taking for example $m = n$ leads to an angle close to $\frac{\pi}{2}$ (since the sampling frequencies are contained in the interval $[-\frac{n}{2} - 2,\frac{n}{2} + 2]$), which necessarily leads to an unstable scenario. Taking more measurements than reconstruction vectors is a common way to stabilize the reconstruction problem
\cite{adcock_2d_wavelet,adcock,ahc1,1,Ma2017,Hrycak,romero,fhgkst,gro,gro1,strom,adnon,wfr,gataricneu,gataricneuneu,MR1882684}.
The operator norm of $Q_\lambda$ is decreasing in
$\lambda $, and the approximation becomes less sensitive to noise,
with  $Q_1$ being the most stable reconstruction in line with  \cite[Theorem
6.2.]{adhapo12} and with  Theorem \ref{simplet}. The intermediate
reconstruction operators offer a  trade-off between sensitivity to
noise and  the  out-of-space contributions. A suitable choice of
$\lambda $ then 
 leads to more accurate reconstructions than $Q_0$ and $Q_1$. For example for $m=10$ and $\ac{SNR} =
20$dB the average relative reconstruction error is $0.088$ for $Q_0$,
$0.098$ for $Q_1$, but only $0.075$ for $Q_{0.1}$.  

Interestingly, for  $\ac{SNR} = 10$dB  we obtain a higher
average approximation error for dimension $m=40$ than for   $m =
10,20,30$. Although the increase in dimension makes the distance
$\|f-P_\cT f\|_{L^2}$ smaller,  the high irregularity in the sampling frequencies
seems to  lead to unstable scenarios. This correlates  nicely with
the  increase of the operator norm $\|Q_\lambda\|_{\rm op}$ with increasing $m$.  

We next  observe that  the condition number
$\kappa(\Sigma_\lambda^{-\frac{1}{2}} U^*T)$  is increasing with
$\lambda $, as
anticipated in  Lemma \ref{cond}. As discussed in Section
\ref{calc}, this opens the possibility of finding non-diagonal
weight matrices. 

\begin{figure*}
	\centering
	\begin{subfigure}[t]{0.45\textwidth}
		\includegraphics[width=\textwidth]{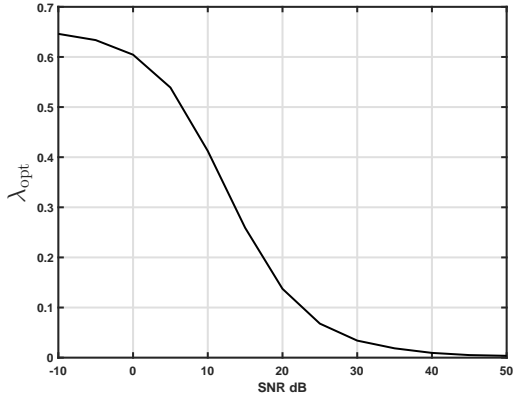}
		\caption{The \ac{SNR} versus $\lambda_{\rm{opt}}$. }
	\end{subfigure}
	\hfill
	\begin{subfigure}[t]{0.45\textwidth}
		\includegraphics[width=\textwidth]{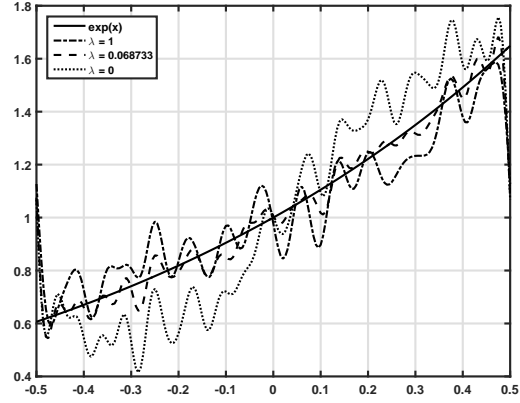}
		\caption{The approximations
obtained by $Q_0$, $Q_{\lambda_{\rm opt}}$ and $Q_1$ have a relative error of $0.151$, $0.059$ and $0.100$
respectively; $\mathrm{SNR} = 20$dB.}
	\end{subfigure}
\begin{subfigure}[t]{0.45\textwidth}
		\includegraphics[width=\textwidth]{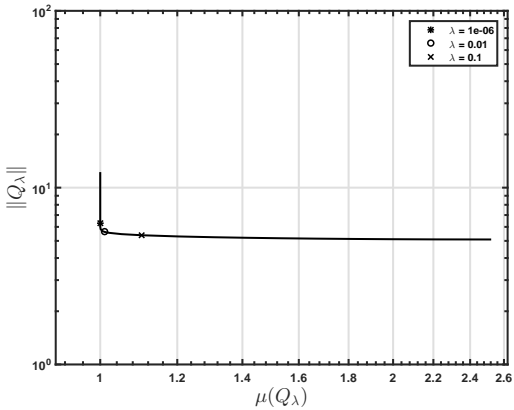}
		\caption{The quasi-optimality constant
                  $\mu(Q_\lambda)$ versus the operator norm
                  $\|Q_\lambda\|_{\rm op}$ for $\lambda \in [0,1]$.} 
	\end{subfigure}
\hfill
	\begin{subfigure}[t]{0.45\textwidth}
		\includegraphics[width=\textwidth]{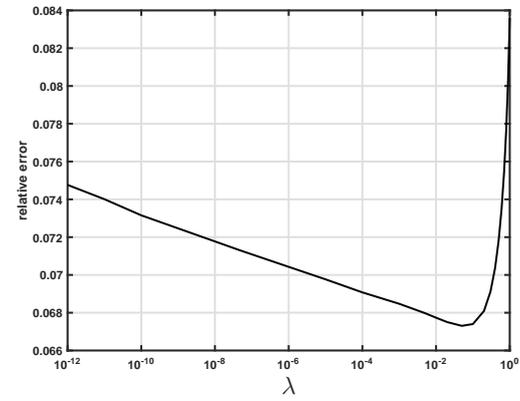}
		\caption{$\lambda$ versus the relative reconstruction error $\frac{\|Q_\lambda d-f\|_{L^2}}{\|f\|_{L^2}}$; $\mathrm{SNR} = 20$dB}
	\end{subfigure}
	\caption{Reconstruction performance of the operators $Q_\lambda$.} \label{fig1}
\end{figure*}

For  Figure \ref{fig1}(a) we have determined the parameter $\lambda_
{\rm opt}$ such that the relative reconstruction error  is minimal, i.e,
$\lambda_
{\rm opt} = \mathrm{arg~min}  \frac{\|Q_\lambda d-f\|_{L^2}}{\|f\|_{L^2}}$. We
then plot the correlation between $\lambda_
{\rm opt}$ and the signal-to-noise
ratio. The plot confirms Theorems~\ref{mini} and~\ref{simplet}: for $\mathrm{SNR}
\to \infty $ the optimal reconstruction is with $Q_0$, and for
$\mathrm{SNR} \to -\infty $, the optimal reconstruction is with $Q_1$.  

In Figure \ref{fig1}(b) we depict the approximations obtained by
$Q_0$,  $Q_1$ and $Q_{\lambda_{\rm{opt}}}$ for a single  realization of the
sampling frequencies and noise with $m = 20$ and \ac{SNR} $= 20$dB.
The optimal choice of the regularization parameter $\lambda_{\rm{opt}}$ yields 
a significantly better approximation than the standard least square
approximation (with $Q_1$). 


In Figure \ref{fig1}(c) we depict the quasi-optimality constant
$\mu(Q_\lambda)$ versus the operator norm $\|Q_\lambda\|_{\rm op}$ for $\lambda
\in [0,1]$  for  $m = 20$ and \ac{SNR}$ = 20$dB. The curve consists of
the points $(\mu (Q_\lambda ), \|Q_\lambda \|_{\rm op})$ for $\lambda \in
[0,1]$. This curve exhibits a  striking change of direction at a small
value of $\lambda $, as is shown by the points corresponding to
$\lambda = 10^{-6}, 0.01, 0.1$. The value of 
$\lambda$ near  the edge of the curve (around $0.01$) yields a good
trade-off between quasi-optimality and operator norm. 

 Figure \ref{fig1}(d) shows the relative reconstruction error
 $\frac{\|Q_\lambda d-f\|_{L^2}}{\|f\|_{L^2}}$ as a function of $\lambda $ for
 fixed $\mathrm{SNR}$. This
 is a typical $L$-curve known from many regularization procedures of
 ill-posed problems. The plots supports the  interpretation
 of $\lambda $ as a regularization parameter.  
\subsection{Reconstruction from measurements of a biased function}\label{befo}
We now assume that we are given a set of Fourier measurements of a perturbation of $f\in \H$
\begin{equation*}
\tilde{d} = [ \mathcal{F}(f+\Delta f)(\omega_{-n}),\dots,\mathcal{F}(f+\Delta f)(\omega_{n})]^T.
\end{equation*} 
The sampling frequencies $\omega_j$ are as in \eqref{omegaj} with $\delta_j\in[-2,2]$.
For each set of sampling frequencies we choose $\Delta f$ as a trigonometric polynomial
\begin{equation}\label{noisy_function}
\Delta f = \sum_{j = -\frac{n}{2}}^{\frac{n}{2}} a_j \e^{2 i \pi j \cdot}\chi_{[-1/2,1/2]}.
\end{equation}
The coefficients $a_j$ in \eqref{noisy_function}
are i.i.d. Gaussian distributed with variance (average power) $\sigma ^2$. 
Table \ref{table2} shows the relative error
$\frac{\|Q_\lambda \tilde{d}-f\|_{L^2}}{\|f\|_{L^2}}$ for $\text{SNR} = \infty$,
$\text{SNR} = 20$dB and $\text{SNR} = 10$dB 
and $m= 10$ in (a), $m = 20$ in (b), $m = 30$ in (c) and $m = 40$ in
(d). In this case the most  accurate reconstruction is always given
by the reconstruction with the  operator
$Q_0$, thus confirming Corollary \ref{tbvc}. In addition, with
increasing dimension of the reconstruction space the distance 
$\|f-P_\cT f\|_{L^2}$ is decreasing with  $m$, and the relative error
decreases up to $m=40$.   
For $m = 50$ the angle between the sampling and reconstruction space
is $\varphi_{\mathcal{T},\mathcal{U}}=0.309 \pm 0.094$, which results
in a significantly higher approximation error. 

\begin{table}[h!]
\centering
\begin{subtable}[t]{0.5\textwidth}
\tiny
\centering
\begin{tabular}{|c | c c  c |}
\toprule
 &\multicolumn{3}{c|}{relative error $\frac{\|Q_\lambda h-f\|_{L^2}}{\|f\|_{L^2}}$}  \\ 
 & $\text{SNR}\! = \!\infty$ & $\text{SNR}\! = \! 20$dB   & $\text{SNR}\! = \! 10$dB  \\
\midrule
		$Q_{0}$	  &$0.067  \pm    0.000 $&$   0.075   \pm   0.003  $&$   0.127  \pm     0.020$  \\
		$Q_{0.1}$	&$0.068  \pm    0.003 $&$   0.076    \pm  0.004  $&$   0.128   \pm    0.020$  \\
		$Q_{0.2}$	&$0.069  \pm    0.007 $&$   0.077    \pm  0.007  $&$   0.129   \pm     0.021$  \\
		$Q_{0.3}$	&$0.071  \pm    0.011 $&$   0.079    \pm   0.011  $&$   0.131  \pm     0.022$  \\
		$Q_{0.4}$	&$0.073  \pm    0.015 $&$   0.081    \pm   0.015  $&$   0.132  \pm     0.024$  \\
		$Q_{0.5}$	&$0.075  \pm    0.019 $&$   0.083    \pm   0.018  $&$   0.134  \pm     0.027$  \\
		$Q_{0.6}$	&$0.078  \pm    0.023 $&$   0.086    \pm   0.022  $&$   0.136  \pm     0.030$  \\
		$Q_{0.7}$	&$0.081  \pm    0.027 $&$   0.088    \pm   0.027  $&$   0.139  \pm      0.034$  \\
		$Q_{0.8}$	&$0.084  \pm    0.032 $&$   0.092    \pm   0.031  $&$   0.142  \pm     0.038$  \\
		$Q_{0.9}$	&$0.088  \pm    0.037 $&$   0.095    \pm   0.037  $&$   0.146  \pm     0.044$  \\
		$Q_{1}$	  &$0.093  \pm    0.043 $&$    0.100   \pm    0.043 $&$   0.151  \pm     0.051$  \\
\bottomrule
\end{tabular}
\caption{$m = 10$}
\end{subtable}%
\begin{subtable}[t]{0.5\linewidth}
\tiny
\centering
\begin{tabular}{|c | c c  c |}
\toprule
 &\multicolumn{3}{c|}{relative error $\frac{\|Q_\lambda h-f\|_{L^2}}{\|f\|_{L^2}}$}  \\ 
 & $\text{SNR}\! = \!\infty$ & $\text{SNR}\! = \! 20$dB    & $\text{SNR}\! = \! 10$dB  \\
\midrule
		$Q_{0}$	  &$0.048     \pm       0.000  $&$      0.067    \pm     0.005    $&$     0.159   \pm      0.023$  \\
		$Q_{0.1}$	&$0.049     \pm  0.002  $&$      0.068     \pm    0.006    $&$      0.160    \pm     0.022$  \\
		$Q_{0.2}$	&$0.050      \pm 0.005   $&$     0.069      \pm   0.007     $&$    0.161     \pm    0.023$  \\
		$Q_{0.3}$	&$0.052     \pm   0.008  $&$      0.071    \pm     0.009    $&$     0.162   \pm      0.023$  \\
		$Q_{0.4}$	&$0.054     \pm   0.011  $&$      0.073    \pm     0.011    $&$     0.163   \pm      0.024$  \\
		$Q_{0.5}$	&$0.056     \pm   0.014  $&$      0.074    \pm     0.013    $&$     0.164   \pm      0.025$  \\
		$Q_{0.6}$	&$0.059     \pm   0.017  $&$      0.076    \pm     0.016    $&$     0.166   \pm      0.027$  \\
		$Q_{0.7}$	&$0.061     \pm    0.020  $&$      0.079    \pm     0.019    $&$     0.168   \pm      0.029$  \\
		$Q_{0.8}$	&$0.064     \pm   0.023  $&$      0.081    \pm     0.022    $&$      0.170   \pm      0.032$  \\
		$Q_{0.9}$	&$0.067     \pm   0.026  $&$      0.084    \pm     0.026    $&$     0.173   \pm      0.037$  \\
		$Q_{1}$	  &$0.071     \pm    0.030  $&$      0.088    \pm      0.030    $&$     0.176   \pm      0.043$  \\
\bottomrule
\end{tabular}
\caption{$m = 20$}

\end{subtable}
\begin{subtable}[t]{0.5\textwidth}
\tiny
\centering
\begin{tabular}{|c | c c  c |}
\toprule
 &\multicolumn{3}{c|}{relative error $\frac{\|Q_\lambda h-f\|_{L^2}}{\|f\|_{L^2}}$}  \\ 
 & $\text{SNR}\! = \!\infty$ & $\text{SNR}\! = \! 20$dB    & $\text{SNR}\! = \! 10$dB  \\
\midrule
		$Q_{0}$	  &$0.039    \pm    0.000   $&$     0.070     \pm    0.006    $&$     0.189    \pm     0.023$  \\
		$Q_{0.1}$	&$0.041    \pm    0.003   $&$     0.071     \pm    0.007    $&$     0.189    \pm     0.023$  \\
		$Q_{0.2}$	&$0.042    \pm    0.006   $&$     0.072     \pm    0.008    $&$     0.190    \pm     0.023$  \\
		$Q_{0.3}$	&$0.044    \pm    0.008   $&$     0.073     \pm    0.009    $&$     0.191    \pm     0.023$  \\
		$Q_{0.4}$	&$0.046    \pm    0.011   $&$     0.075     \pm    0.010    $&$     0.192    \pm     0.024$  \\
		$Q_{0.5}$	&$0.048    \pm    0.013   $&$     0.076     \pm    0.012    $&$     0.193    \pm     0.024$  \\
		$Q_{0.6}$	&$0.050    \pm    0.015   $&$     0.078     \pm    0.014    $&$     0.194    \pm     0.025$  \\
		$Q_{0.7}$	&$0.052    \pm    0.017   $&$     0.079     \pm    0.016    $&$     0.195    \pm     0.026$  \\
		$Q_{0.8}$	&$0.055    \pm    0.020   $&$     0.081     \pm    0.018    $&$     0.197    \pm     0.027$  \\
		$Q_{0.9}$	&$0.057    \pm    0.022   $&$     0.083     \pm    0.020    $&$     0.199    \pm     0.029$  \\
		$Q_{1}$	  &$0.060    \pm    0.025   $&$     0.086     \pm    0.023    $&$     0.202    \pm     0.033$  \\
\bottomrule
\end{tabular}
\caption{$m = 30$}
\end{subtable}%
\begin{subtable}[t]{0.5\linewidth}
\tiny
\centering
\begin{tabular}{|c | c c  c |}
\toprule
 &\multicolumn{3}{c|}{relative error $\frac{\|Q_\lambda h-f\|_{L^2}}{\|f\|_{L^2}}$}  \\ 
 & $\text{SNR}\! = \!\infty$ & $\text{SNR}\! = \! 20$dB    & $\text{SNR}\! = \! 10$dB  \\
\midrule
		$Q_{0}$	  &$0.034   \pm     0.000   $&$     0.075    \pm     0.007   $&$      0.215    \pm     0.023$  \\
		$Q_{0.1}$	&$0.036   \pm     0.003   $&$     0.076    \pm     0.007   $&$      0.216    \pm     0.023$  \\
		$Q_{0.2}$	&$0.038   \pm     0.006   $&$     0.077    \pm     0.008   $&$      0.217    \pm     0.023$  \\
		$Q_{0.3}$	&$0.040   \pm     0.008   $&$     0.078    \pm     0.008   $&$      0.218    \pm     0.023$  \\
		$Q_{0.4}$	&$0.041   \pm     0.010   $&$     0.079    \pm     0.009   $&$      0.218    \pm     0.023$  \\
		$Q_{0.5}$	&$0.043   \pm     0.011   $&$     0.080    \pm     0.010   $&$      0.219    \pm     0.024$  \\
		$Q_{0.6}$	&$0.045   \pm     0.013   $&$     0.082    \pm     0.011   $&$      0.220    \pm     0.024$  \\
		$Q_{0.7}$	&$0.046   \pm     0.014   $&$     0.083    \pm     0.013   $&$      0.221    \pm     0.025$  \\
		$Q_{0.8}$	&$0.048   \pm     0.016   $&$     0.084    \pm     0.014   $&$      0.222    \pm     0.026$  \\
		$Q_{0.9}$	&$0.050   \pm     0.017   $&$     0.085    \pm     0.015   $&$      0.223    \pm     0.027$  \\
		$Q_{1}$	  &$0.051   \pm     0.019   $&$     0.087    \pm     0.016   $&$      0.225    \pm     0.029$  \\
\bottomrule
\end{tabular}

\caption{$m = 40$}

\end{subtable}
		\caption{Reconstruction of the exponential function by
                  the 
                  operators $Q_\lambda$ from measurements of the
                  biased function }\label{table2} 
\end{table}

\section*{\bf\uppercase{Appendix: Frames in Hilbert Spaces}}

We need the definition of the Moore-Penrose pseudoinverse of an
operator on a Hilbert space \cite[Section 2.5]{ch08}.
We use the notation $\R(A)$ for the range, and $\N(A)$ for the null-space of the operator $A$.
\begin{defi}\label{pseudoinverse}
Let $\H$ and $\W$ be Hilbert spaces. If $A:\W\rightarrow \H$ is a
bounded operator with a closed range $\R(A)$, then there exists a  unique bounded operator $A^\dagger:\H\rightarrow \W$ satisfying 
\begin{align*}
 &\N(A^\dagger) = \R(A)^\perp = \N(A^*), \\ 
 &\R(A^\dagger) = \N(A)^\perp = \R(A^*), \text{ and }\\ 
 &A A^\dagger x =x,~x\in \R(A).
 \end{align*}
The operator $A^\dagger$ is called the Moore-Penrose pseudoinverse of $A$.
\end{defi}

For a sequence  $(u_j)_{j\in \mathbb{N}}$ we define the synthesis
operator on the subspace of finite sequences by 
$$
U(c_j)_{j\in \mathbb{N}}=\sum_{j=1}^\infty c_j u_j \, .
$$

\begin{defi} \label{frfr}
(i)  If $U$ can be extended to a bounded operator $   U:\ell
^2(\mathbb{N})\rightarrow \H$, $(u_j)_{j\in \bN } $ is called a Bessel
sequence.

(ii) If $U$ is bounded $   U:\ell
^2(\mathbb{N})\rightarrow \H$, $(u_j)_{j\in \bN } $ and has closed
range, $(u_j)_{j\in \bN } $ is called a frame for the subspace $\cU =
\overline{\textnormal{span}}(u_j)_{j\in \mathbb{N}}$.  

(iii) If $(u_j)_{j\in \bN } $ is a Bessel sequence, then the adjoint
operator of $U$ is the \textit{analysis operator}
\begin{equation*}
   U^*:\H\rightarrow \ell ^2(\mathbb{N}),\quad  U^*f= (\langle
   f,u_j\rangle_\cH)_{j\in \mathbb{N}} \, , 
\end{equation*}
and $S = UU^*:\cH \to \cH $, $Sf = \sum_{j=1}^\infty \langle
f,u_j\rangle_\cH u_j$ is the frame operator of $(u_j)$. 

(iv) The sequence $(S^{\dagger} u_j)_{j\in \bN }  \subseteq
\cU $ is the canonical dual frame in $\cU$ , and every $f\in \cU $  possesses
the frame expansions
 $$f = \sum_{j\in \mathbb{N}} \langle f,S^{\dagger}u_j\rangle_\cH u_j =
 \sum_{j\in \mathbb{N}} \langle f,u_j\rangle_\cH S^{\dagger}u_j
$$ 
with  unconditional convergence of both series. 
\end{defi}

\begin{lem} \label{tight}
Let $\mathcal{U}$ be a closed subspace of $\H$ and let
$(u_j)_{j\in\mathbb{N}}$ be a frame for $\mathcal{U}$.
The set 
\begin{equation*}
 (S^{\frac{\dagger}{2}} u_j)_{j\in\mathbb{N}}
\end{equation*}
forms a tight frame for $\mathcal{U}$ with frame bound equal to $1$. 
The synthesis operator $M$ of the  sequence $(S^{\frac{\dagger}{2}} u_j)_{j\in\mathbb{N}}$
is given by $M := S^{\frac{\dagger}{2}} U$, and
\begin{equation*} 
 P_\mathcal{U} = MM^* = S^{\frac{\dagger}{2}}S S^{\frac{\dagger}{2}} = S^{\dagger}S = S S^{\dagger}.
\end{equation*}    
\end{lem}

Lemma \ref{tight1} proves the following. Suppose that we are given the inner products $(\langle f,u_j\rangle_\cH)_{j\in \mathbb{N}}$ of an element $f\in \H$ with a frame $(u_j)_{j\in \mathbb{N}}$ for $\mathcal{U}$ (a closed subspace of $\H$). Applying the operator $( U^* U)^{\frac{\dagger}{2}}$ to these measurements, we obtain the inner products of $f$ with the tight frame
$(S^{\frac{\dagger}{2}} u_j)_{j\in\mathbb{N}}$ for $\mathcal{U}$.

\begin{lem}\label{tight1}
Let $\mathcal{U}$ be a closed subspace of $\H$ and $(u_j)_{j\in\mathbb{N}}$  a frame for $\mathcal{U}$ with synthesis  operator $ U$, analysis operator $ U^*$, Gramian $G = U^*U$ and frame operator $S = U U^*$. Then
\begin{equation*} 
 G^{\frac{\dagger}{2}} U^* =  U^*S^{\frac{\dagger}{2}}.
 \end{equation*}
Thus, 
 $G^{\frac{\dagger}{2}} U^*$
is the analysis operator of the tight frame sequence $(S^{\frac{\dagger}{2}} u_j)_{j\in\mathbb{N}}$.
\begin{proof}
 Obviously for $k\in \mathbb{N}$
 \begin{equation*}
 ( U^* U)^k U^* =  U^*( U  U^*)^k. 
 \end{equation*}
Therefore, 
\begin{equation*}
 p(G) U^* =  U^*p(S) 
\end{equation*}
for every polynomial $p$. We are going to prove that there exists a sequence of polynomials $(p_k)_{k\in\mathbb{N}}$, such that for $i = 1,2$
\begin{equation*}
 \underset{m\rightarrow \infty} {\lim}\|p_m(M_i) - M_i^\frac{\dagger}{2}\|_{\rm op} = 0
\end{equation*}
simultaneously for $M_1:= G$ and $M_2:= S$. 

Let $A$ and $B$ denote the lower bound and upper frame bound of the frame sequence $(u_j)_{j\in\mathbb{N}}$.
From the lower frame bound $A$ we infer that for every $f\in \mathcal{U} = \N( U U^*)^\perp = \N(S)^\perp$
\begin{equation*}
A \|f\|_\cH^2\leq \langle S f,f\rangle_\cH.                                                          
\end{equation*}
Consequently the set $\sigma(S)\backslash \{0\}$ is bounded below by $A$. Here $\sigma(S)$ denotes the spectrum of the operator $S$. The upper frame bound $B$ ensures that the set $\sigma(S)$ has the upper bound $B$. 
This shows that $0$ is an isolated point of the spectrum, and that for ${K:=\{0\}\cup [A,B]}$ the function 
${h: K\rightarrow \mathbb{R}}$ 
\begin{align*}\
h(x) = 
\begin{cases}
 \frac{1}{\sqrt{x}} &\text{for }x\in [A,B],\\
  0 &\text{for }x= 0
\end{cases} 
\end{align*}
is continuous on $K$. Since $\sigma(S)\cup \{0\} = \sigma(G)\cup \{0\}$, $h$ is also continuous 
on $\sigma(G)$.

By the Weierstrass approximation theorem there exists a sequence of polynomials $(p_m)_{m\in \mathbb{N}}$, such that 
\begin{equation*}
 \underset{m\rightarrow \infty} {\lim}\|p_m -h\|_\infty = 0, 
\end{equation*}
uniformly on $K$.
By the continuous functional calculus
\begin{equation*}
  \underset{m\rightarrow \infty}{\lim} \|p_m(M_i)-h(M_i)\|_{\rm op} = 0
\end{equation*}
simultaneously for $M_1:= G$ and $M_2:= S$ and  $h(M_i) = M_i^\frac{\dagger}{2}$ for $i = 1,2$.
%
%
%
\end{proof}
\end{lem}



\begin{lem}{\cite[Lemma 5.3.6]{ch08}}\label{candual}
Let $(f_j)_{j\in \mathbb{N}}$ be a frame for $\H$ with frame operator $S$ and let $f\in \H$. If $f$ has a representation $f = \sum_{j\in \mathbb{N}} c_j f_j$ for some coefficients $(c_j)_{j\in \mathbb{N}}$, then
\begin{equation*}
\sum_{j\in \mathbb{N}}|c_j|^2 = \sum_{j\in \mathbb{N}}|\langle f,S^{-1}f_j\rangle_\cH|^2 + \sum_{j\in \mathbb{N}}|c_j-\langle f,S^{-1}f_j\rangle_\cH|^2.
\end{equation*}
\end{lem}

\end{document}